\documentclass[11pt]{amsart}
\usepackage{bbm}
\usepackage{mathrsfs}

\usepackage{amscd}
\usepackage{amsmath}
\usepackage{graphicx}
\usepackage{amsfonts}
\usepackage{amssymb}
\begin{document}
\title[Lie ring isomorphisms]
{ Lie ring isomorphisms between nest algebras on Banach spaces}

\author{Xiaofei Qi}
\author{Jinchuan Hou}
\author{Juan Deng}

\address[Qi, Hou and Deng]
{Department of Mathematics, Taiyuan University of Technology,
Taiyuan 030024, P. R. China; Department of Mathematics, Shanxi
University, Taiyuan, 030006, P. R. China.}
\email{xiaofeiqisxu@aliyun.com; jinchuanhou@aliyun.com;}
\email{houjinchuan@tyut.edu.cn; juanhappyforever@163.com}

\thanks{ This work is  supported by National Natural Science
Foundation of China (11171249,11101250, 11271217) and Youth
Foundation of Shanxi Province (2012021004).}

\begin{abstract}
Let ${\mathcal N}$ and ${\mathcal M}$ be nests on
  Banach spaces $X$ and $Y$ over the (real or
complex) field $\mathbb F$ and let $\mbox{\rm Alg}{\mathcal N}$ and
$\mbox{\rm Alg}{\mathcal M}$ be the associated nest algebras,
respectively.  It is shown that a map $\Phi:{\rm Alg}{\mathcal
N}\rightarrow{\rm Alg}{\mathcal M}$ is a Lie ring isomorphism (i.e.,
$\Phi$ is additive, Lie multiplicative and bijective)  if and only
if $\Phi$ has the form $\Phi(A) = TAT^{-1} + h(A)I$ for all $A\in
\mbox{\rm Alg}{\mathcal N}$ or $\Phi(A)=-TA^*T^{-1}+h(A)I$ for all
$A\in \mbox{\rm Alg}{\mathcal N}$, where $h$ is an additive
functional vanishing on all commutators and $T$ is an invertible
bounded linear or conjugate linear operator when $\dim X=\infty$;
$T$ is a bijective $\tau$-linear transformation for some field
automorphism $\tau$ of $\mathbb F$ when $\dim X<\infty$.

\end{abstract}

\maketitle

\section{Introduction and main results}

Let $\mathcal R$ and $\mathcal R^\prime$ be two associative rings.
Recall that a map $\phi: \mathcal R \rightarrow \mathcal R^\prime$
is called a multiplicative map if $\phi(AB)=\phi(A)\phi(B)$ for any
$A,B \in \mathcal R$; is called a Lie multiplicative map if
$\phi([A,B])=[\phi(A),\phi(B)]$ for any $A,B \in \mathcal R$, where
$[A,B]=AB-BA$ is the Lie product of $A$ and $B$ which is also called
a commutator. In addition, a map $\phi: \mathcal R \rightarrow
\mathcal R^\prime$ is called a Lie multiplicative isomorphism if
$\phi$ is bijective and Lie multiplicative;  is called a Lie ring
isomorphism if $\phi$ is bijective, additive and Lie multiplicative.
If $\mathcal R$ and $\mathcal R^\prime$ are algebras  over a field
$\mathbb F$, $\phi:{\mathcal R}\rightarrow{\mathcal R}^\prime $ is
called a Lie algebraic isomorphism if $\phi$ is bijective, $\mathbb
F$-linear and Lie multiplicative. For the study of Lie ring
isomorphisms between rings, see \cite{B3,B5,M1} and the references
therein. In this paper we focus our attention on Lie ring
isomorphisms between nest algebras on general Banach spaces.

Let $X$ be a Banach space over the (real or complex) field ${\mathbb
F}$ with topological dual $X^*$. $\mathcal B(X)$ stands for the
algebra of all bounded linear operators on $X$. A nest $\mathcal N$
on $X$ is a complete totally ordered subspace lattice, that is, a
chain of closed (under norm topology) subspaces of $X$ which is
closed under the formation of arbitrary closed linear span (denote
by $\bigvee$) and intersection (denote by $\bigwedge$), and which
includes $\{0 \}$ and $X$. The nest algebra associated with a nest
$\mathcal N$, denoted by ${\rm Alg}\mathcal N$, is the weakly closed
operator algebra consisting of all operators that leave every
subspace $N\in\mathcal N$ invariant. For $N\in \mathcal N$, let
$N_{-}=\bigvee\{M\in \mathcal N \mid M \subset N\}$ and $
N_{-}^{\perp}=(N_{-})^{\perp}$, where $N^{\perp}=\{f\in X^*\mid N
\subseteq \ker(f)\}$. 
If $\mathcal N$ is a nest on $X$, then ${\mathcal
N}^{\perp}=\{N^\perp\mid N\in{\mathcal N}\}$ is a nest on $X^*$ and
$({\rm Alg}{\mathcal N})^*\subseteq{\rm Alg}{\mathcal N}^\perp$. If
${\mathcal N}=\{(0),X\}$, we say that $\mathcal N$ is a trivial
nest, in this case, Alg${\mathcal N}={\mathcal B}(X)$. Non-trivial
nest algebras are very important reflexive operator algebras that
are not semi-simple, not semi-prime and not self-adjoint. If $\dim
X<\infty$,  a nest algebra on $X$ is isomorphic to an algebra of
upper triangular block matrices. Nest algebras are studied
intensively by a lot of literatures. For more details on basic
theory of nest algebras, the readers can refer to \cite{D,L}.

In \cite{MS}, Marcoux and Sourour proved that every Lie algebraic
isomorphism  between nest algebras   on separable complex Hilbert
spaces is a sum $\alpha+\beta$, where $\alpha$ is an algebraic
isomorphism or the negative of an algebraic anti-isomorphism and
$\beta:{\rm Alg}{\mathcal N}\rightarrow{\mathbb C}I$ is a linear map
vanishing on all commutators, that is,  satisfying $\beta([A, B])=0$
for all $A,B\in{\rm Alg}{\mathcal N}$.

Qi and Hou in \cite{QH} generalized the result of Marcoux and
Sourour by classifying certain Lie multiplicative isomorphisms. Note
that, a Lie multiplicative isomorphism needs not  be   additive. Let
${\mathcal N}$ and ${\mathcal M}$ be nests on Banach spaces $X$ and
$Y$ over the (real or complex) field $\mathbb F$, respectively, with
the property that if $M\in {\mathcal M}$ such that $M_-=M$, then $M$
is complemented in $Y$ (Obviously, this assumption is not needed if
$Y$ is a Hilbert space or if $\dim Y<\infty$). Let $\mbox{\rm
Alg}{\mathcal N}$ and $\mbox{\rm Alg}{\mathcal M}$ be respectively
the associated nest algebras, and let $\Phi:{\rm Alg}{\mathcal
N}\rightarrow{\rm Alg}{\mathcal M}$ be a bijective map. Qi and Hou
in \cite{QH} proved that, if $\dim X=\infty$ and if there is a
nontrivial element in $\mathcal N$ which is complemented in $X$,
then $\Phi$ is a Lie multiplicative isomorphism if and only if there
exists a map $h:\mbox{\rm Alg}{\mathcal N}\rightarrow {\mathbb F}I$
  with $h([A,B])=0$ for all $A, B\in \mbox{\rm Alg}{\mathcal
N}$ such that $\Phi$ has the form $\Phi(A)=TAT^{-1}+h(A) $ for all
$A\in \mbox{\rm Alg}{\mathcal N}$ or $\Phi(A)=-TA^*T^{-1}+h(A) $ for
all $A\in \mbox{\rm Alg}{\mathcal N}$, where, in the first form,
$T:X\rightarrow Y$ is an invertible bounded linear or
conjugate-linear operator so that $N\mapsto T(N)$ is an order
isomorphism from $\mathcal N$ onto $\mathcal M$, while in the second
form, $X$ and $Y$ are reflexive, $T: X^*\rightarrow Y$ is an
invertible bounded linear or conjugate-linear operator so that
$N^\bot\mapsto T(N^\bot)$ is an order  isomorphism from ${\mathcal
N}^\bot$ onto $\mathcal M$. If $\dim X=n<\infty$, identifying nest
algebras with upper triangular block matrix algebras, then $\Phi$ is
a Lie multiplicative isomorphism if and only if there exist a field
automorphism $\tau : {\mathbb F}\rightarrow{\mathbb F}$ and certain
invertible matrix $T$ such that either
$\Phi(A)=TA_{\tau}T^{-1}+h(A)$ for all $A$, or
$\Phi(A)=-T(A_{\tau})^{\rm tr}T^{-1}+h(A)$ for all $A$, where
$A_\tau =(\tau (a_{ij}))$ for $A=(a_{ij})$ and $A^{\rm tr}$ is the
transpose of $A$. Particularly, above results give
  a characterization of Lie ring isomorphisms between nest algebras
  for finite-dimensional case, and for infinite-dimensional case
  under the mentioned assumptions on $\mathcal N$ and $\mathcal M$.

Recently, Wang and Lu in \cite{WL} generalized Marcoux and Sourour's
result from another direction, and proved that every Lie algebraic
isomorphism between nest algebras ${\rm Alg}{\mathcal N}$ and ${\rm
Alg}{\mathcal M}$ for any nests  $\mathcal N$ and $\mathcal M$   on
Banach spaces $X$ and $Y$ respectively can be decomposed as
$\alpha+\beta$, where $\alpha$ is an algebraic isomorphism or the
negative of an algebraic anti-isomorphism and $\beta:{\rm
Alg}{\mathcal N}\rightarrow{\mathbb F}I$ is a linear map vanishing
on each commutator. Because  Lie algebraic isomorphisms were
characterized  in \cite{BE} for the case that the nest $\mathcal N$
has a nontrivial  complemented element, Wang and Lu  in \cite{WL}
mainly dealt with the case that all nontrivial elements of $\mathcal
N$ are not complemented.

\if Thus, it is interesting and natural to ask how to characterize
and classify all Lie ring isomorphisms between nest algebras of
Banach space operators? Based on the methods in \cite{QH} and
\cite{WL}, \fi

The purpose of the present paper is to
  characterize  all Lie ring isomorphisms
between   nest algebras of Banach space operators for any nests.
Note that, the Lie ring isomorphisms are very different from
algebraic ones. For example, the method used in \cite{BE} to
characterize Lie algebraic isomorphisms for the case that the nest
$\mathcal N$ has a nontrivial complemented element is not valid for
characterizing Lie ring isomorphisms. Algebraic isomorphisms between
nest algebras are continuous, however ring isomorphisms are not
necessarily continuous for finite-dimensional case \cite[Remark
2.6]{HZ}.

The following are the main results of this paper.

{\textbf  {Theorem 1.1.} {\it Let ${\mathcal N}$ and ${\mathcal M}$
be nests on  Banach spaces $X$ and $Y$ over the (real or complex)
field $\mathbb F$, and, ${\rm Alg}{\mathcal N}$ and ${\rm
Alg}{\mathcal M}$ be the associated nest algebras, respectively.
Then a map $\Phi:{\rm Alg}{\mathcal N}\rightarrow{\rm Alg}{\mathcal
M}$ is a Lie ring isomorphism, that is, $\Phi$ is additive,
bijective and satisfies $\Phi([A,B])=[\Phi(A),\Phi(B)]$ for all
$A,B\in{\rm Alg}{\mathcal N}$, if and only if $\Phi$ has the form
$\Phi(A)=\Psi(A)+h(A)I$ for all $A\in{\rm Alg}{\mathcal N}$, where
$\Psi$ is a ring isomorphism or the negative of a ring
anti-isomorphism between the nest algebras and $h:\mbox{\rm
Alg}{\mathcal N}\rightarrow {\mathbb F}$ is an additive functional
satisfying $h([A,B])=0$ for all $A, B\in \mbox{\rm Alg}{\mathcal
N}$.}

The ring isomorphisms and the   ring anti-isomorphisms between nest
algebras of Banach space operators were characterized in
\cite[Theorem 2.2, Theorem 2.7 and Remark 2.6]{HZ}. Using these
results and Theorem 1.1, we can get more concrete characterization
of Lie ring isomorphisms. Recall that a map $S:W\rightarrow V$ with
$W,V$ linear spaces over a field $\mathbb F$ is called $\tau$-linear
if $S$ is additive and $S(\lambda x)=\tau(\lambda)Sx$ for all $x\in
W$ and $\lambda\in{\mathbb F}$, where $\tau$ is a field automorphism
of $\mathbb F$.

{\bf Theorem 1.2.} {\it Let ${\mathcal N}$ and ${\mathcal M}$ be
nests on Banach spaces $X$ and $Y$ over the (real or complex) field
$\mathbb F$,  and let ${\rm Alg}{\mathcal N}$ and ${\rm
Alg}{\mathcal M}$ be the associated nest algebras, respectively.
Then a map $\Phi:{\rm Alg}{\mathcal N}\rightarrow{\rm Alg}{\mathcal
M}$   is a Lie ring isomorphism if and only if there exist an
additive functional $h:\mbox{\rm Alg}{\mathcal N}\rightarrow
{\mathbb F} $ satisfying $h([A,B])=0$ for all $A, B\in \mbox{\rm
Alg}{\mathcal N}$ and a field automorphism $\tau: {\mathbb
F}\rightarrow{\mathbb F}$ such that one of the following holds.}

(1) {\it There exists a $\tau$-linear transformation $T: X
\rightarrow Y$ such that the map $N\mapsto T(N) $ is an order
isomorphism from $\mathcal N$ onto $\mathcal M$ and
$$\Phi(A)=TAT^{-1}+h(A)I \ \ {for \ all} \ \ A \in {\mbox{\rm Alg}} \mathcal
N.$$ }

(2) {\it $X$ and $Y$ are reflexive, there exists a $\tau$-linear
transformation  $T: X^* \rightarrow Y$ such that the map $N_-
^{\perp}\mapsto T(N_- ^{\perp}) $ is an order isomorphism from $
{\mathcal N }^ {\perp}$ onto ${ \mathcal M}$ and
$$\Phi(A)=-TA^*T^{-1}+h(A)I\ \ { for \ all}\ \ A \in {\mbox{\rm Alg}} \mathcal
N.$$}

{\it Moreover, if $\dim X=\infty$, the above $T$ is in fact an
invertible bounded linear or conjugate-linear operator; if ${\mathbb
F}={\mathbb R}$, $T$ is linear.}

For the finite dimensional case, it is clear that every nest algebra
on a finite dimensional space is isomorphic to an upper triangular
block matrix algebra. Let ${\mathcal M}_n({\mathbb F})$ denote the
algebra of all $n\times n$ matrices over $\mathbb{F}$. Recall that
an upper triangular block matrix algebra ${\mathcal T}={\mathcal
T}(n_1, n_2,\ldots ,n_k)$ is a subalgebra of $ {\mathcal
M}_n({\mathbb F})$ consisting of all $n\times n$ matrices of the
form
$$A=\left(\begin{array} {cccc}
A_{11}& A_{12}&\ldots&A_{1k}\\
0& A_{22}&\ldots& A_{2k}\\
\vdots&\vdots&\ddots&\vdots\\
0&0&\ldots& A_{kk}
\end{array}\right),$$ where $n_1, n_2, \cdot\cdot\cdot, n_k$ are finite sequence of positive
integers satisfying $n_1+n_2+ \cdot\cdot\cdot+ n_k=n$ and
$A_{ij}\in{\mathcal M}_{n_i\times n_j}({\mathbb F})$, the space of
all $n_i\times n_j$ matrices over $\mathbb F$. Thus by Theorem 1.2,
we get a characterization of Lie ring isomorphisms between upper
triangular block matrix algebras.

{\bf Corollary 1.3.} {\it Let ${\mathbb F}$ be the real or complex
field, and $m, n$ be positive integers greater than 1. Let
${\mathcal T}={\mathcal T}(n_1, n_2,\ldots ,n_k)\subseteq {\mathcal
M}_n({\mathbb F})$ and ${\mathcal S}={\mathcal T}(m_1, m_2,\ldots
,m_r)\subseteq {\mathcal M}_m({\mathbb F})$ be upper triangular
block matrix algebras, and $\Phi:{\mathcal T}\rightarrow{\mathcal
S}$ be a
 map. Then $\Phi$ is a Lie ring isomorphism if and only if
$m=n$, and there exist an additive functional $\phi:{\mathcal
T}\rightarrow {\mathbb F}$ satisfying $\phi([A,B])=0$ for all $A,
B\in {\mathcal T}$, a field automorphism $\tau:{\mathbb
F}\rightarrow{\mathbb F}$ such that either}

(1) {\it ${\mathcal T}={\mathcal S}$, there exists   an invertible
matrix $T\in {\mathcal T}$ such that
$$\Phi(A)=TA_{\tau}T^{-1}+\phi(A)I \ \  for \ \ all \ \ A\in{\mathcal T};$$
or}

(2)  {\it $(n_1, n_2,\ldots ,n_k)=(m_r, m_{r-1},\ldots ,m_1)$, there
exists   an invertible block matrix $T=(T_{ij})_{k\times k}$ with
$T_{ij}\in {\mathcal M}_{n_i\times n_j}({\mathbb F})$ and $T_{ij}=0$
whenever $i+j> k+1$, such that
$$\Phi(A)=-TA_{\tau}^{\rm tr}T^{-1}+\phi(A)I \ \  for\ \  all \ \ A\in{\mathcal T}.$$
Where $A_\tau=(\tau(a_{ij}))_{n\times n}$ for $A=(a_{ij})_{n\times
n}\in {\mathcal M}_n({\mathbb F})$ and $A^{\rm tr}$ is the transpose
of $A$. If ${\mathbb F}={\mathbb R}$, then $\Phi$ is a Lie algebraic
isomorphism.}

Corollary 1.3 is also a consequence of \cite[Corollary 2.2]{QH}.

\if false Theorem 1.1 allows us give a generalization of the main
result obtained in \cite{QH}, too.

 {\bf Theorem  1.3.} {\it Let ${\mathcal N}$ and ${\mathcal M}$ be nests on  Banach spaces $X$
and $Y$ over the (real or complex) field $\mathbb F$, and, ${\rm
Alg}{\mathcal N}$ and ${\rm Alg}{\mathcal M}$ be the associated nest
algebras, respectively. Assume that $\mathcal N$ contains a
complemented non-trivial element. Then a map $\Phi:{\rm
Alg}{\mathcal N}\rightarrow{\rm Alg}{\mathcal M}$ is a Lie
multiplicative isomorphism, that is, $\Phi$ is   bijective and
satisfies $\Phi([A,B])=[\Phi(A),\Phi(B)]$ for all $A,B\in{\rm
Alg}{\mathcal N}$, if and only if $\Phi$ has the form $\Phi=\Psi+h$,
where $\Psi$ is a ring isomorphism or the negative of a ring
anti-isomorphism between the nest algebras and $h:\mbox{\rm
Alg}{\mathcal N}\rightarrow {\mathbb F}I$ is a  map satisfying
$h([A,B])=0$ for all $A, B\in \mbox{\rm Alg}{\mathcal N}$.}

Thus, similar to Corollary 1.2, by the characterization of ring
isomorphisms and ring anti-isomorphisms between nest algebras on
Banach spaces (for both finite-dimensional and infinite-dimensional
cases, ref. \cite{HZ}), one can get a classification of all Lie
multiplicative isomorphisms between nest algebras on Banach spaces
under the assumption that the nests contain a non-trivial
complemented element.\fi

Since the Lie ring isomorphisms between nest algebras on
finite-dimensional Banach spaces were already characterized in
\cite[Corollary 2.2]{QH}, to give a   classification of all Lie ring
isomorphisms between nest algebras of Banach space operators, it
suffices to prove Theorem 1.1 for  the infinite-dimensional cases
without any additional assumption on the nests. It is clear  that
$\dim X=\infty\Leftrightarrow\dim Y=\infty$.

The remaining part of the paper is to prove the main result Theorem
1.1 under the assumption that both $X$ and $Y$ are
infinite-dimensional. Our approach borrows and  combines some ideas
developed in \cite{QH} and \cite{WL}. In Section 2, we give
preliminary lemmas and some of them are also parts of the proof of
the main result. Section 3 deals with the case that both $(0)$ and
$X$ are limit points of the nest $\mathcal N$, that is, $(0)=(0)_+$
and $X_-=X$. The case that $X_-\not= X$ and $X_-$ is complemented or
$(0)\not= (0)_+$ and $(0)_+$ is complemented is discussed in Section
4. And finally, the case that $X_-\not= X$ and $X_-$ is not
complemented or $(0)\not= (0)_+$ and $(0)_+$ is not complemented is
considered  in Section 5.

\section{Preliminary lemmas}

In this section, we give some preliminary lemmas, definitions and
symbols which are needed in other sections to prove the main
result.

Let $X$ and $Y$ be Banach spaces over $\mathbb F$, and let $\mathcal
N$ and $\mathcal M$ be nests on $X$ and $Y$. Let Alg$\mathcal N$ and
Alg$\mathcal M$ be associated nest algebras, respectively. It is
well known that the commutant of a nest algebra is trivial, i.e., if
$T\in{\mathcal B}(X)$ and $TA=AT$ for every operator $A
\in$Alg$\mathcal N$, then $T=\lambda I$ for some scalar
$\lambda\in\mathbb F.$ This fact will be used in this paper without
any specific explanation. In addition, the symbols ${\rm ran} T$,
$\ker T$ and ${\rm rank}T$ stand for the range, the kernel and the
rank (i.e., the dimension of ran$T$) of an operator $T$,
respectively. For $x\in X$ and $f\in X^*$, $x\otimes f$ stands for
the operator on $X$ with rank not greater than 1 defined by
$(x\otimes f)y=f(y)x$ for every $y$. Some times we use $\langle
x,f\rangle$ to present the value $f(x)$ of $f$ at $x$.

The following lemma is also well known which gives a
characterization of rank one operators in nest algebras.

{\bf Lemma 2.1.} (\cite{L}) {\it  Let ${\mathcal N}$ be a nest on a
(real or complex) Banach space $X$ and $x\in X$, $f\in X^*$. Then
 $x\otimes f\in{\rm Alg}\mathcal N$ if and only if there exists a
subspace $N\in{\mathcal N}$ such that $x\in N$ and $f\in
N_-^\perp$.}

For any non-trivial element $E\in\mathcal{N}$, define
$$\mathcal{J}({\mathcal N},E)=\{ A\in {\rm Alg}{\mathcal N}: AE=0\
\ {\rm and}\ \ A^*E^{\perp}=0\}.\eqno(2.1)$$ In \cite{WL}, Wang and
Lu proved that $\mathcal L$  is a proper maximal commutative Lie
algebra ideal in ${\rm Alg}\mathcal N$  if and only if  ${\mathcal
L}=\mathbb{F} I +{\mathcal J}({\mathcal N},E)$ for some unique
$E\in{\mathcal N}$. The following lemma shows that any maximal
commutative Lie ring ideal also  rises  in this way.

{\textbf {Lemma 2.2.} {\it ${\mathcal J}$ is a proper maximal
commutative Lie ring ideal in ${\rm Alg}{\mathcal N}$ if and only if
it is   a proper maximal commutative Lie algebra ideal.}

{\bf Proof.} Assume that $\mathcal{J}$ is a maximal commutative Lie
ring ideal. Then for any $A\in{\rm Alg}{\mathcal N}$, any
$C\in{\mathcal J}$ and any $\lambda\in \mathbb{F}$, we have
$[A,\lambda C]=\lambda [A,C]\in \mathbb{F} \mathcal{J}$, which
implies that $\mathbb{F} \mathcal{J}$ is a Lie ring ideal. It is
obvious that $\mathbb{F} \mathcal{J}$ is also commutative. So
$\mathbb{F}\mathcal{J}\subseteq \mathcal{J}$ as $\mathcal J$ is
maximal. Note that $\mathbb{F} {\mathcal J}\supseteq {\mathcal J}$.
Thus we get $\mathbb{F} \mathcal{J} = \mathcal{J}$. It follows that
$\mathcal{J}$ is also a Lie algebra ideal. The converse is obvious.
\hfill$\Box$

In the rest part of this paper, we  assume that $\Phi:{\rm
Alg}{\mathcal N}\rightarrow{\rm Alg}{\mathcal M}$ is a Lie ring
isomorphism. If ${\mathcal N}$ contains at least one nontrivial
element,   by Lemma 2.2, for any nontrivial element
$E\in\mathcal{N}$, $\Phi(\mathbb{F} I +\mathcal{J}(\mathcal{N},E))$
is a maximal commutative Lie ring ideal in ${\rm Alg}\mathcal M$.
Hence there is a unique nontrivial element $F\in  {\mathcal M}$ such
that $\Phi(\mathbb{F} I +\mathcal{J}(\mathcal{N},E))=\mathbb{F} I
+\mathcal{J}(\mathcal{M},F).$ Define a map
$$\hat{\Phi}:{\mathcal N}\setminus{\{(0),X\}}\rightarrow {\mathcal M}\setminus{\{(0),Y\}} \eqno(2.2)$$ by $\Phi(\mathbb{F} I
+\mathcal{J}(\mathcal{N},E))=\mathbb{F} I
+\mathcal{J}(\mathcal{M},\hat{\Phi}(E))$.

With the symbols introduced above and by an argument similar to
\cite[Lemmas 4.1, 4.3, 4.4]{WL}, one can show that the following
lemma is still true for the Lie ring isomorphism $\Phi$.

{\bf Lemma 2.3.}  {\it $\hat{\Phi}$ in Eq.(2.2) is bijective and is
either order-preserving or order-reversing, that is, $\hat{\Phi}$ is
an order isomorphism or a reverse-order isomorphism from ${\mathcal
N}$ onto ${\mathcal M}$ if we extend the definition of $\hat{\Phi}$
so that $\hat{\Phi}((0))=(0)$ or $Y$ and $\hat{\Phi}(X)=Y$ or $(0)$
accordingly. }

By Lemma 2.3, for any $A\in{\mathcal J}({\mathcal N},E)$ with
nontrivial $E\in{\mathcal N}$, there exists a unique operator
$B\in{\mathcal J}({\mathcal M},\hat{\Phi}(E))$ such that
$\Phi(A)-B\in{\mathbb F}I$. Thus we can define another map
$$\bar{\Phi}:\bigcup\{{\mathcal J}({\mathcal N},E): E\in{\mathcal N}\ {\rm is\ nontrivial}\}\rightarrow
\bigcup\{{\mathcal J}({\mathcal M},F): F\in{\mathcal M}\ {\rm is\
nontrivial}\} \eqno(2.3)$$ with the property that
$\bar{\Phi}(A)\in{\mathcal J}({\mathcal M},\hat{\Phi}(E))$ and
$\Phi(A)-\bar{\Phi}(A)\in{\mathbb F}I$ for any $A\in{\mathcal
J}({\mathcal N},E)$.

Similar to \cite[Lemma 4.2]{WL}, we have

{\bf Lemma 2.4.}  {\it  $\bar{\Phi}$ is a bijective map and
$\bar{\Phi}({\mathcal J}({\mathcal N},E))={\mathcal J}({\mathcal
M},\hat{\Phi}(E))$ for every non-trivial $E\in{\mathcal N}$.}

Next we discuss the idempotents in nest algebras. Denote by
$\mathcal{E}(\mathcal N)$  the set  of all idempotents in ${\rm
Alg}{\mathcal N}$.

{\textbf {Lemma 2.5.} (\cite[Lemma 2.2]{QH}) {\it Let ${\mathcal N}$
be a nest  on a (real or complex) Banach space $X$  and $A\in{\rm
Alg}\mathcal N$.}

(1) {\it $A \in {\mathbb F}I+\mathcal{E}({\mathcal N})$  if and only
if $[A,[A,[A,T]]]=[A,T]$  for all $T \in {\rm Alg}\mathcal N$.}

(2) {\it $A$ is the sum of a scalar and an idempotent operator with
range in $\mathcal N$   if and only if $[A,[A,T]]=[A,T]$  for all
$T\in {\rm Alg}\mathcal N$.}

By Lemma 2.5,  if $P$ is an idempotent operator in ${\rm
Alg}{\mathcal N}$, then $\Phi(P)=Q+\lambda_P I$, where $\lambda_P\in
{\mathbb F}$ and $Q$ is an idempotent operator in ${\rm
Alg}{\mathcal M}$. Furthermore, if ${\rm ran}P\in {\mathcal N}$,
then ${\rm ran}Q\in{\mathcal M}$. So we can define a map
$$\tilde{\Phi}: \mathcal{E}(\mathcal N)\rightarrow \mathcal{E}(\mathcal M) \eqno(2.4)$$ by
$\tilde{\Phi}(P)=\Phi(P)-\lambda_P I$. It is easily seen that
$\tilde{\Phi}$ is a bijective map from ${\mathcal E}({\mathcal N})$
onto ${\mathcal E}({\mathcal M})$; see \cite{QH}.

Now, for any nontrivial element $E\in {\mathcal N}$,  define two
sets
$$\Omega_{1}(\mathcal{N},E)=\{P\in \mathcal{E}(\mathcal{N}):
PE=0\}\quad{\rm and}\quad \Omega_{2}(\mathcal{N},E)=\{P\in
\mathcal{E}(\mathcal{N}): P^{*}E^{\bot}=0\}. \eqno(2.5)$$ For any
nontrivial element $F\in{\mathcal M}$, the sets
$\Omega_{1}(\mathcal{M},F)$ and $\Omega_{2}(\mathcal{M},F)$ can be
analogously defined. Note that, if $P\in\Omega_1({\mathcal N},E)$,
then one can easily check $P^*E^\perp\not=0$; if $E$ is not
complemented, then $PE=0\Rightarrow(I-P)^*E^\perp\not=0$ and
$P^*E^\perp=0\Rightarrow(I-P)E\not=0$. These facts are needed in the
proof of Lemma 2.6.

Still, by an  argument similar to \cite[Lemmas 5.2-5.5]{WL}, one can
show that the following lemma is  true for the Lie ring isomorphism
$\Phi$, with $\hat{\Phi}$ and $\tilde{\Phi}$ defined in Eq.(2.2) and
Eq.(2.4) respectively.

{\bf Lemma 2.6.} {\it Assume that $E\in{\mathcal N}$ is  nontrivial
and   not complemented in $X$. }

(1) {\it Either $\tilde{\Phi}(\Omega_1({\mathcal
N},E))\subseteq\Omega_1({\mathcal N},\hat{\Phi}(E))$ or
$I-\tilde{\Phi}(\Omega_1({\mathcal N},E))\subseteq\Omega_2({\mathcal
N},\hat{\Phi}(E))$, and either $\tilde{\Phi}(\Omega_2({\mathcal
N},E))\subseteq\Omega_2({\mathcal N},\hat{\Phi}(E))$ or
$I-\tilde{\Phi}(\Omega_2({\mathcal N},E))\subseteq\Omega_1({\mathcal
N},\hat{\Phi}(E))$.}

(2) {\it If $\Omega_1({\mathcal N},\hat{\Phi}(E))$ and
$\Omega_2({\mathcal N},\hat{\Phi}(E))$ are not empty, then
$\tilde{\Phi}(\Omega_1({\mathcal N},E))\subseteq\Omega_1({\mathcal
N},\hat{\Phi}(E))$ if and only if  $\tilde{\Phi}(\Omega_2({\mathcal
N},E))\subseteq\Omega_2({\mathcal N},\hat{\Phi}(E))$.}

(3) {\it If $F\in{\mathcal N}$ is also nontrivial  and  not
complemented in $X$ with $F<E$, then   $\Omega_1({\mathcal
N},E)\not=\emptyset$ and $\tilde{\Phi}(\Omega_1({\mathcal
N},E))\subseteq\Omega_1({\mathcal N},\hat{\Phi}(E))$ together imply
that $\tilde{\Phi}(\Omega_1({\mathcal
N},F))\subseteq\Omega_1({\mathcal N},\hat{\Phi}(F))$. }

The following lemma gives a characterization of complemented
elements $E\in{\mathcal N}$ by the operators in ${\mathcal
J}({\mathcal N},E)$ and ${\mathcal E}({\mathcal N})$, which is
needed to prove that $\hat{\Phi}$ preserves the complementarity.

{\bf Lemma 2.7.} {\it Assume that $E\in{\mathcal N}$ is a nontrivial
element. The following statements are equivalent.}

(1) {\it $E\in{\mathcal N}$ is complemented in $X$.}

(2) {\it There exists some idempotent $P\in{\rm Alg}{\mathcal N}$
 such that $[P,A]=A$ for any $A\in{\mathcal J}({\mathcal N},E)$.}

(3) {\it There exists some idempotent $P\in{\rm Alg}{\mathcal N}$
 such that $[P,A]=A$ for any rank one operator $A\in{\mathcal J}({\mathcal N},E)$.}

{\it Furthermore,   we have $E={\rm ran}P$ for $P$ in (2) and (3).}

{\bf Proof.} (1)$\Rightarrow$(2). If $E\in{\mathcal N}$ is
complemented in $X$, there exists an idempotent $P\in {\rm
Alg}{\mathcal N}$ such that ${\rm ran}P=E$. \if According to the
space decomposition $X=E+\ker P$, for any $A\in{\mathcal
J}({\mathcal N},E)$, we have
$$P=\left(\begin{array}{cc} I& 0\\ 0 & 0\end{array}\right)\quad{\rm and}\quad
A=\left(\begin{array}{cc} 0& A_{12}\\ 0 & 0\end{array}\right).$$\fi
For any $A\in{\mathcal J}({\mathcal N},E)$, we have $PA=A$ and
$AP=0$. Hence $[P,A]=A$.

(2)$\Rightarrow$(3) is obvious.

(3)$\Rightarrow$(1). Assume that  there exists some idempotent
$P\in{\rm Alg}{\mathcal N}$ such that $[P,A]=A$ for any rank-1
operator $A\in{\mathcal J}({\mathcal N},E)$. According to the space
decomposition $X={\rm ran}P+\ker P$, for any $A\in{\mathcal
J}({\mathcal N},E)$, we have
$$P=\left(\begin{array}{cc} I& 0\\ 0 & 0\end{array}\right)\quad{\rm and}\quad
A=\left(\begin{array}{cc} A_{11}& A_{12}\\ A_{21} &
A_{22}\end{array}\right).$$ Then $[P,A]=A$ implies that
$A=\left(\begin{array}{cc} 0& A_{12}\\ 0 & 0\end{array}\right)$. So
$PA=A$ and $AP=0$ hold for all $A\in{\mathcal J}({\mathcal N},E)$.
Take any $y\in E$ and any $g\in E^\perp$. It is easy to check that
$y\otimes g\in{\mathcal J}({\mathcal N},E)$. So $Py\otimes
g=y\otimes g$ and $y\otimes gP=0$.  It follows that $Py=y$ and
$P^*g=0$ for all $y\in E$ and all $g\in E^\perp$. Thus we get
$E\subseteq{\rm ran} P$ and $E^\bot\subseteq\ker P^*$. If
$E\not={\rm ran} P$, then there exist $x\in{\rm ran}P$ and $g\in
E^\bot$ such that $\langle x,g\rangle=1$. This leads to a
contradiction that $0=\langle x, P^*g\rangle=\langle Px,
g\rangle=\langle x,  g\rangle=1$. Hence we must have ${\rm ran}P=E$
and $E$ is complemented in $X$. \hfill$\Box$

{\bf Lemma 2.8.} {\it Non-trivial element $E\in{\mathcal N}$ is
complemented in $X$ with ${\rm ran}P=E$ if and only if
$\hat{\Phi}(E)$ is complemented in $Y$ with ${\rm
ran}\tilde{\Phi}(P)=\hat{\Phi}(E)$. Here $P\in{\rm Alg}{\mathcal
N}$ is an idempotent.}

{\bf Proof.} It is clear that
$\widehat{\Phi^{-1}}=\hat{\Phi}^{-1}$. So we  need only to  show
that $\hat{\Phi}(E)$ is complemented in $Y$ whenever
$E\in{\mathcal N}$ is complemented in $X$. Indeed, if
$E\in{\mathcal N}$ is complemented, by Lemma 2.7, there exists
some idempotent $P\in{\rm Alg}{\mathcal N}$ such that $[P,A]=A$
for any $A\in{\mathcal J}({\mathcal N},E)$. By definitions of
$\hat{\Phi}$,
 $\bar{\Phi}$ and $\tilde{\Phi}$ (ref.  Eqs.(2.2)-(2.4)), there exists some scalar
$\gamma$ such that
$$\bar{\Phi}(A)+\gamma
I=\Phi(A)=\Phi([P,A])=[\Phi(P),\Phi(A)]
=[\tilde{\Phi}(P),\bar{\Phi}(A)].$$ Since $\bar{\Phi}$ maps
${\mathcal J}({\mathcal N},E)$ onto ${\mathcal J}({\mathcal
M},\hat{\Phi}(E))$ by Lemma 2.4, we see that $B+\gamma
I=[\tilde{\Phi}(P),B] $ holds for any $B\in{\mathcal J}({\mathcal
M},\hat{\Phi}(E))$.  Assume that $B$ is of rank-1. If
$\gamma\not=0$, then  $[\tilde{\Phi}(P),B]=B+\gamma I$ is a sum of
nonzero scalar and a rank-1 operator, which is impossible since a
commutator can not be the sum of a nonzero scalar and a compact
operator. Hence $\gamma=0$, and $B=[\tilde{\Phi}(P),B]$ for all
rank-1 operators $B\in{\mathcal J}({\mathcal M},\hat{\Phi}(E))$. By
Lemma 2.7 again, $\hat{\Phi}(E)$ is complemented in $Y$ with ${\rm
ran}\tilde{\Phi}(P)=\hat{\Phi}(E)$. \hfill$\Box$

The following lemma is obvious.

{\bf Lemma 2.9.} {\it $\Phi({\mathbb F}I)={\mathbb F}I$.}

Finally, we give a lemma, which is needed to prove our main result.
Let $E$ and $F$ be subspaces of $X$ and $X^*$, respectively. Denote
by $E\otimes F$ the set $\{x\otimes f:x\in E,f\in F\}$.

{\bf Lemma  2.10.} {\it Let $X_{i}$ be an infinite-dimensional
Banach space, $i=1,2$. Let $E_{i}$ and $F_{i}$ be closed subspaces
with dimensions $> 2$ of $X_{i}$ and $X_{i}^*$ respectively. Let
$\mathcal{A}_{i}$ be a unital subalgebra of ${\mathcal B}(X_{i})$
containing $E_{i}\otimes F_{i}$. Suppose that $\Psi:{\mathcal
A}_1\rightarrow{\mathcal A}_2$ is an additive bijective map
satisfying $\Psi(\mathbb F I)=\mathbb F I$ and  $\Psi({\mathbb F}
I+E_{1}\otimes F_{1})={\mathbb F} I+E_{2}\otimes F_{2}$. Then there
is a map $\gamma: E_{1}\otimes F_{1} \rightarrow \mathbb{F}$ and a
field automorphism $\tau: \mathbb{F} \rightarrow \mathbb{F}$ such
that either}

(1) {\it $\Psi(x\otimes f) = \gamma(x,f) I + Cx\otimes Df$ for all
$x\in E_{1}$ and $f\in F_{1},$ where $C: E_{1} \rightarrow E_{2}$
and $D: F_{1} \rightarrow F_{2}$ are two $\tau$-linear bijective
maps; or}

(2) {\it  $\Psi(x\otimes f)=\gamma(x,f) I+Df\otimes Cx$ for all
$x\in E_{1}$ and $f\in F_{1},$ where  $C: E_{1} \rightarrow F_{2}$
and $D: F_{1} \rightarrow E_{2}$ are two $\tau$-linear bijective
maps.}

Lemma 2.10 can be proved by a similar approach as that in \cite{BHD}
and we omit its proof here.

Since the ``if'' part of Theorem 1.1 is obvious, we  need only to
check the ``only if'' part. If $\mathcal N$ is a trivial nest, that
is, ${\mathcal N}=\{(0),X\}$, then ${\mathcal M}=\{(0),Y\}$ by Lemma
2.3. So $\Phi$ is a Lie ring isomorphism from ${\mathcal B}(X)$ onto
${\mathcal B}(Y) $. Bai, Du and Hou showed in \cite{BDH} that every
Lie multiplicative isomorphism between prime rings with a
non-trivial idempotent element is of the form $\psi+\beta$ with
$\psi$ a ring isomorphism or the negative of a ring anti-isomorphism
and $\beta$ a central valued map vanishing on each commutator. Note
that ${\mathcal B}(X)$ is prime and contains non-trivial idempotents
if $\dim X\geq 2$. Hence, for the case that $\mathcal N$ is trivial,
Theorem 1.1 follows from \cite{BDH}.

So, in the rest sections  we always assume that ${\mathcal N}$ is
nontrivial. Thus ${\mathcal M}$ is also nontrivial by Lemma 2.3. We
shall complete the proof of Theorem 1.1 by considering several cases
according to the situations of $(0)_+$ and $X_-$, which will be
dealt with separably in Sections 3-5.

\section{The case that $(0)_{+}=(0)$ and $X_{-}=X$}

In this section, we deal with the case that both $(0)$ and $X$ are
limit points of $\mathcal N$, i.e., $(0)_{+}=(0)$ and $X_{-}=X$.

Keep the definitions of $\hat{\Phi}$, $\bar{\Phi}$ in mind; ref.
Eqs.(2.2) and (2.3).

{\bf Lemma 3.1.} {\it Assume that $(0)_{+}=(0)$ and $X_{-}=X$. Let
$E\in{\mathcal N}$ be a nontrivial element. Then for any $x\in E$
and $f\in E^{\bot}$, $\bar{\Phi}(x\otimes f)$ is a rank one
operator.}

{\bf Proof.} Since, by Lemma 2.3, the bijective map
$\hat{\Phi}:{\mathcal N}\rightarrow{\mathcal M}$ is either
order-preserving or order-reversing, we have  $(0)_+=(0)$ and
$Y_-=Y$ in $\mathcal{M}$. Write $R=\bar{\Phi}(x\otimes f)$ with
$\bar{\Phi}$ defined in Eq.(2.3). Then $R\in
\mathcal{J}(\mathcal{M},\hat{\Phi}(E))$. If, on the contrary, ${\rm
rank}R\geq 2$, then, since $\bigcup\{M: M\in{\mathcal M}\}$ is a
dense linear manifold of $Y$, there exists a nontrivial element
$M_0\in{\mathcal M}$ and two vectors $u,v\in M_0$ such that $Ru$ and
$Rv$ are linearly independent. As $\bigcap\{M:(0)<M\in{\mathcal
M}\}=\{0\}$, there exists some nontrivial  $L\in{\mathcal M}$ such
that $Ru,Rv\not\in L$. Let $Y_L={\rm span}\{Ru,L\}$. By Hahn-Banach
theorem, there exists some $g\in Y_L^\bot$ such that $g(Rv)\not=0$.
Then $g\in L^{\bot}$ satisfying $g(Ru)=0$ and $g(Rv)\neq0.$ It is
easily checked that $L<\hat{\Phi}(E)< M_0$. Take $z\in L$, $h\in
M_0^{\bot}$, and  let $A=\bar{\Phi}^{-1}(z\otimes g)$,
$B=\bar{\Phi}^{-1}(u\otimes h)$. By Lemma 2.4, $A\in
\mathcal{J}(\mathcal{N},\hat{\Phi}^{-1}(L))$ and $B\in
\mathcal{J}(\mathcal{N},\hat{\Phi}^{-1}(M_0)).$

If $\hat{\Phi}$ is order-preserving, we have
$\hat{\Phi}^{-1}(L)<E<\hat{\Phi}^{-1}(M_0)$. So
$$\Phi(A(x\otimes f)B)=\Phi([A,[x\otimes f,B]])=[z\otimes g,[R,u\otimes h]]=(z\otimes g)R(u\otimes
h)=0.$$ It follows from the injectivity of $\Phi$ that $A(x\otimes
f)B=0$, which implies $A(x\otimes f)=0$ or $(x\otimes f)B=0.$ If
$A(x\otimes f)=0,$ then $0=\Phi([A,x\otimes f])=[z\otimes
g,R]=(z\otimes g)R\not=0$, a contradiction; if $(x\otimes f)B=0,$
then $0=\Phi([x\otimes f,B])=[R,u\otimes h]=R(u\otimes h)\not=0,$ a
contradiction.

If $\hat{\Phi}$ is  order-reversing, we have
$\hat{\Phi}^{-1}(M_0)<E<\hat{\Phi}^{-1}(L).$  So$$\Phi(B(x\otimes
f)A)=\Phi([B,[x\otimes f,A]])=[u\otimes h,[R,z\otimes g]]=(z\otimes
g)R(u\otimes h)=0,$$ which yields that either $B(x\otimes f)=0$ or
$(x\otimes f)A=0$. By a similar argument to that of the above, one
can obtain a contradiction.

Therefore, we must have ${\rm rank}\bar{\Phi}(x\otimes f)={\rm
rank}R=1$. \hfill$\Box$ \vspace{3mm}

{\bf Proof of Theorem 1.1 for the case that $(0)=(0)_+$ and
$X_-=X$.}

By  Lemma 3.1, Lemma 2.9 and the bijectivity of $\bar{\Phi}$, we
have proved that, for every nontrivial element $E\in{\mathcal N}$,
  $\Phi(\mathbb{F}I + E\otimes E^{\bot})=\mathbb{F} I +
\hat{\Phi}(E)\otimes \hat{\Phi}(E)^{\bot}$. So, by Lemma 2.10, there
exists  a ring automorphism $\tau_E: \mathbb{F} \rightarrow
\mathbb{F}$ and a map $\gamma_E: E\otimes E^{\bot}\rightarrow
\mathbb{F}$ such that either
$$\Phi(x\otimes f) =
\gamma_E(x,f) I + C_Ex\otimes D_{E^\perp}f \eqno (3.1)$$ holds for
all $x\in E$ and $f\in E^{\bot}$, where   $C_E: E\rightarrow
\hat{\Phi}(E)$ and $D_{E^\perp}: E^{\bot} \rightarrow
\hat{\Phi}(E)^{\bot}$ are two $\tau$-linear bijective maps; or
$$\Phi(x\otimes f)=\gamma_E(x,f)I +D_{E^\perp}f\otimes C_Ex \eqno (3.2)$$
holds for all $x\in E$ and $f\in E^{\bot}$, where  $C_E: E
\rightarrow \hat{\Phi}(E)^\perp$ and
 $D_{E^\perp}:E^{\bot} \rightarrow  \hat{\Phi}(E)$ are two $\tau$-linear bijective
maps.

It is easily checked that, if there is a nontrivial $E_0\in{\mathcal
N}$ such that Eq.(3.1) holds, then  Eq.(3.1) holds for any
nontrivial $E\in{\mathcal N}$; If there is a nontrivial
$E_0\in{\mathcal N}$ such that Eq.(3.2) holds, then  Eq.(3.1) holds
for any nontrivial $E\in{\mathcal N}$.

Assume that Eq.(3.1) holds for a nontrivial $E\in{\mathcal N}$.
Then, for any $N\in{\mathcal N}$, any $x\in E\cap N$ and any $f\in
N^{\perp}\cap E^{\perp}$, we have
$$\Phi(x\otimes f)=\gamma_E(x,f)I + C_Ex\otimes D_{E^\perp}f=\gamma_N(x,f)I + C_Nx\otimes D_{N^\perp}f.$$
Since ${\rm rank}I=\infty$, the above equation yields
$\gamma_N(x,f)=\gamma_E(x,f)$ and $C_Ex\otimes
D_{E^\perp}f=C_Nx\otimes D_{N^\perp}f$ for any $x\in N\cap E$. This
entails that there exists an automorphism $\tau: {\mathbb
F}\rightarrow{\mathbb F}$ so that $\tau_N=\tau_E=\tau$ for any
$N,E\in{\mathcal N}$, and if $N\subset E$, then there exists a
scalar $\alpha_{EN}$ such that $C_E|_N=\alpha_{EN}C_N$ and
$D_{N^{\perp}}|_{E^{\perp}}=\alpha_{EN}D_{E^{\perp}}$. Now fix
$E\in{\mathcal N}$.
 For any $N\in{\mathcal N}$, we define
$$\left\{
\begin{array}{lll}
\tilde{C}_N=C_N, & \tilde{D}_{N_{-}^{\perp}}=D_{N^{\perp}}, & \mbox{\rm if}\ \ N=E; \\
\tilde{C}_N=\alpha_{EN}C_N, &
 \tilde{D}_{N_{-}^{\perp}}=\frac{1}{\alpha_{EN}}D_{N^{\perp}}, & \mbox{\rm if} \ \ N\subset E; \\
\tilde{C}_N=\frac{1}{\alpha_{EN}}C_N, &
 \tilde{D}_{N_{-}^{\perp}}=\alpha_{EN}D_{N^{\perp}}, & \mbox{\rm if} \ \
N\supset E.
\end{array}
\right.$$ It is obvious that $\{\tilde{C_N}: N\not=(0),X\}$ and $\{
\tilde{D}_{N_{-}^{\perp}}: N\not=(0),X\}$ are well defined with
$\tilde{C}_N|_E=\tilde{C}_E$ and
$\tilde{D}_{E^{\perp}}|_{N^\perp}=\tilde{D}_{N^{\perp}}$ whenever
$E\subseteq N$. Thus there exist bijective $\tau$-linear maps $C:
\bigcup\{N\in{\mathcal N}: N\not=(0),X\}\rightarrow
\bigcup\{M\in{\mathcal M}: M\not=(0),Y\} $ and $D: \bigcup\{N^\perp:
N\in{\mathcal N}, N\not=(0),X\}\rightarrow\bigcup\{M^\perp:
M\in{\mathcal M}, M\not=(0),Y\}$ such that $C| _{N}=\tilde{C}_N$ and
$D|_{N^{\perp}}=\tilde{D}_{N_-^{\perp}}$ for any $N\in
{\mathcal{N}}\backslash\{(0),X\}$.

By now, we have shown that, there exist bijective $\tau$-linear maps
$C: \bigcup\{N\in{\mathcal N}: N\not=(0),X\}\rightarrow
\bigcup\{M\in{\mathcal M}: M\not=(0),Y\} $ and $D: \bigcup\{N^\perp:
N\in{\mathcal N}, N\not=(0),X\}\rightarrow\bigcup\{M^\perp:
M\in{\mathcal M}, M\not=(0),Y\}$, and a map $\gamma:\bigcup\{E\times
E^\perp: E\in{\mathcal N}\setminus\{(0),X\}\}\rightarrow{\mathbb
F}$, such that for any $x\in  N $ and $f\in N^\perp$ with
$N\in{\mathcal N}\setminus\{(0),X\}$, we have
$$\Phi(x\otimes f)=\gamma(x,f) I+Cx\otimes Df.\eqno(3.3)$$
Therefore, for   any $A\in {\rm Alg}\mathcal{N}$, any $x\in N $ and
any  $f\in N^\perp$, by Eq.(3.3), we have
$$ \Phi([A,x\otimes f])=[\Phi(A),\Phi(x\otimes f)]=\Phi(A)Cx\otimes Df - Cx\otimes
\Phi(A)^{*}Df$$and
 $$\begin{array}{rl}\Phi([A,x\otimes f])=&\Phi(Ax\otimes f-x\otimes
A^*f)\\=&(\gamma(Ax, f)- \gamma(x, A^{*}f))I + CAx\otimes Df -
Cx\otimes DA^{*}f.\end{array}$$ Combining the above two equations
and noting that $I$ is of infinite-rank, one obtains that
$$Cx\otimes \Phi(A)^{*}Df-Cx\otimes DA^{*}f = \Phi(A)Cx\otimes Df -CAx\otimes Df$$
holds for any $x\in N$, $f\in N^\bot$ and any nontrivial
$N\in{\mathcal N}$. Note that $D$ is bijective. So there exists a
scalar $h(A)$ such that
$$\Phi(A)Cx=CAx+h(A)Cx \eqno(3.4)$$ for all $x\in \bigcup\{N\in{\mathcal N}:
N\not=(0),X\}$. It is clear that $h$ is additive as a functional of
Alg$\mathcal N$. Define $\Psi(A)=\Phi(A)-h(A)I$ for all $A\in{\rm
Alg}{\mathcal N}$. Then, by Eq.(3.4), for any $A,B\in{\rm
Alg}{\mathcal N}$ and any $x\in\bigcup\{N\in{\mathcal N}:
N\not=(0),X\}$, we have
$$\Psi(AB)Cx=CABx=\Psi(A)CBx=\Psi(A)\Psi(B)Cx.$$
Since $\bigcup\{N\in{\mathcal N}: N\not=(0),X\}$ is   dense in $X$
and $C$ is bijective, it follows that $\Psi(AB)= \Psi(A)\Psi(B)$ for
all $A,B\in {\rm Alg}\mathcal{N}$, that is, $\Psi$ is a ring
isomorphism and $\Phi(A)=\Psi(A)+h(A)I$ for all $A\in{\rm
Alg}\mathcal N$.

Similarly, if Eq.(3.2) holds, one can check that $\Phi(A)=-\Psi(A) +
h(A)I$ for all $A\in {\rm Alg}\mathcal{N}$, where $\Psi$ is a ring
anti-isomorphism and $h$ is an additive functional.

This completes the proof of Theorem 1.1 for the case that
$(0)_+=(0)$ and $X_-=X$.\hfill$\Box$

\section{The case  $X_-\not=X$ and $X_-$ is complemented or $(0)_+\not=0$ and $(0)_+$ is complemented}

We give only the proof in detail for the case that $X_-\not=X$ and
$X_-$ is complemented. The case that $(0)\not=(0)_+ $ and $(0)_+$ is
complemented in $X$ can be dealt with similarly. \vspace{3mm}

{\bf Proof of Theorem 1.1 for the case that  $  X_{-}\not=X$ and
$X_{-}$ is complemented in $X$.}

Assume that  $  X_{-}\not=X$ and $X_{-}$ is complemented in $X$.

Since $X_{-}$ is complemented, there exists an idempotent $P_0\in
{\rm Alg} \mathcal{N}$ such that ${\rm ran}P_0 = X_{-}$. \if The
proof is similar to that of \cite[Theorem 3.1]{QH}. For
completeness, we give the sketch of the proof here.\fi

With $\tilde{\Phi}$ as in Eq.(2.4) and by a similar argument to that
in the proof of \cite[Theorem 3.1]{QH}, one can show that there
exist an idempotent operator $Q_0$ and a scalar $\lambda_{P_0}$ such
that $\Phi(P_0)=Q_0+\lambda_{P_0} I$ with $Q_0=\tilde{\Phi}(P_0)$,
${\rm ran}Q_0\in \mathcal{M}$ and the following statements hold:

(a)  If there is an idempotent $P_1\in {\rm Alg}{\mathcal N}$ such
that $P_1< P_0$ and $\tilde{\Phi}(P_1)< Q_0$ (or $P_1>P_0$ and
$\tilde{\Phi}(P_1)> Q_0$), then for any $P\in {\rm Alg}{\mathcal
N}$, $P<P_0\Rightarrow\tilde{\Phi}(P)< Q_0$ and
$P>P_0\Rightarrow\tilde{\Phi}(P)>Q_0$.

(b) If there is an idempotent $P_1\in {\rm Alg}{\mathcal N}$ such
that $P_1<P_0$ and $\tilde{\Phi}(P_1)> Q_0$ (or $P_1>P_0$ and
$\tilde{\Phi}(P_1)<Q_0$), then for any $P\in {\rm Alg}{\mathcal N}$,
$P<P_0\Rightarrow\tilde{\Phi}(P )>Q_0$ and
$P>P_0\Rightarrow\tilde{\Phi}(P )< Q_0$.}

{\bf Claim 4.1.} If  (a) occurs,  then $\Phi=\Psi+ h$, where
$\Psi:{\rm Alg}{\mathcal N}\rightarrow {\rm Alg}{\mathcal M}$ is a
ring isomorphism and $h:{\rm Alg}{\mathcal N}\rightarrow{\mathbb
F}I$ is an additive map vanishing on all commutators.

(a) implies that ${\rm ran}Q_0=Y_-$ by Lemmas 2.3 and 2.8. For the
convenience, let ${\mathcal A}_{11}=P_0({\rm Alg}{\mathcal N})P_0$,
${\mathcal A}_{12}=P_0({\rm Alg}{\mathcal N})(I-P_0)$, ${\mathcal
A}_{22}=(I-P_0)({\rm Alg}{\mathcal N})(I-P_0)$, ${\mathcal
B}_{11}=Q_0({\rm Alg}{\mathcal M})Q_0$, ${\mathcal B}_{12}=Q_0({\rm
Alg}{\mathcal M})(I-Q_0)$ and ${\mathcal B}_{22}=(I-Q_0)({\rm
Alg}{\mathcal M})(I-Q_0)$. Then ${\rm Alg}{\mathcal N}={\mathcal
A}_{11}\dot{+}{\mathcal A}_{12}\dot{+}{\mathcal A}_{22}$ and ${\rm
Alg}{\mathcal M}={\mathcal B}_{11}\dot{+}{\mathcal
B}_{12}\dot{+}{\mathcal B}_{22}$, where $\dot{+}$ stands for the
algebraic direct sum.

We will prove Claim 4.1 by several steps.

{\bf Step 1.} {\it $\Phi({\mathcal A}_{12})={\mathcal B}_{12}$. }

The proof is the same as that of \cite[Lemma 2.8]{QH}.

{\bf Step 2.} {\it $\Phi({\mathcal A}_{ii})\subseteq{\mathcal
B}_{ii}+{\mathbb F}I$, $i=1,2$. }

For any $A_{11}\in{\mathcal A}_{11}$, denote
$\Phi(A_{11})=S_{11}+S_{12}+S_{22}$, where $S_{ij}\in {\mathcal
B}_{ij}$. Then
$$0=\Phi([A_{11}, P_0])
=[\Phi(A_{11}),\Phi(P_0)]=[\Phi(A_{11}), Q_0],$$ which implies that
$S_{12}=0$. Let $P\in {\mathcal A}_{22}$ be any idempotent with
$P\not=I-P_0$. It is clear that $P<(I-P_0)$. Then $I-P>P_0$. As
$\Phi$ meets  (a), we have $\tilde{\Phi}(I-P)>Q_0$, that is,
$\tilde{\Phi}(P)<I-Q_0$. It follows from $[A_{11}, P]=0$ that
$[\Phi(A_{11}),\tilde{\Phi}(P)]=[S_{22},\tilde{\Phi}(P)]=0$. By the
arbitrariness of $P$  and the bijectivity of $\tilde{\Phi}$, we see
that $S_{22}$ commutes with every idempotent in ${\mathcal B}_{22}$.
Note that ${\rm ran}Q_0=Y_-$. So $\mathcal{B}_{22} =
\mathcal{B}(\ker Q_0)$, which implies $S_{22}\in {\mathbb F}
(I-Q_0)$. It follows that $\Phi(A_{11})=S_{11}+\lambda
(I-Q_0)=(S_{11}-\lambda Q_0)+\lambda I$ and hence,
$\Phi(\mathcal{A}_{11})\subseteq \mathcal{B}_{11} + \mathbb{F} I$.

Assume that $A_{22}\in \mathcal{A}_{22}$. In the same way as above,
one can show that, $\Phi(A_{22})=T_{11}+T_{22}$ for some
$T_{ii}\in{\mathcal B}_{ii}$, $i=1,2$, with $T_{11}$ commuting with
every idempotent in ${\mathcal B}_{11}$. For any $B_{11}\in
\mathcal{B}_{11}$, by the surjectivity of $\Phi$, there exists some
$A_{0}\in {\rm Alg}\mathcal{N}$ such that $\Phi(A_{0})=B_{11}$.
Furthermore, $A_{0}=A_{11}+\lambda I$ for some $A_{11}\in{\mathcal
A}_{11}$ and some scalar $\lambda$   by Step 1. Thus we have
$$[B_{11},T_{11}] =[B_{11},T_{22}+T_{11}]=[\Phi(A_{0}), \Phi(A_{22})] = \Phi([A_{0},A_{22}]) = \Phi([A_{11}
+\lambda I,A_{22}])=0$$ for all $B_{11}\in \mathcal{B}_{11}$, which
implies $T_{11}\in \mathbb{F} I$ as ${\mathcal A}_{11}$ is a nest
algebra. So $\Phi(\mathcal{A}_{22})\subseteq \mathcal{B}_{22} +
\mathbb{F} I$.

By Steps 1-2, for each $A_{12}\in {\mathcal A}_{12}$, there exists
$B_{12}\in {\mathcal B}_{12}$ such that $\Phi(A_{12})=B_{12}$; for
each $A_{ii}\in {\mathcal A}_{ii}$, $i=1,2$, there exist $B_{ii}\in
{\mathcal B}_{ii}$ and $\lambda_{ii} \in {\mathbb F}$ such that
$\Phi(A_{ii})=B_{ii}+\lambda_{ii} I$.  We claim that $B_{ii}$ and
$\lambda_{ii}$ are uniquely determined. In fact, if
$\Phi(A_{ii})=B_{ii}+\lambda_{ii} I=B_{ii}^{'}+\lambda^{'}_{ii} I$,
then $B_{ii}-B_{ii}^{'}\in {\mathbb F}I$, which implies that
$B_{ii}=B_{ii}^{'}$ and $\lambda_{ii}=\lambda^{'}_{ii}$. Let
$\Psi(A_{ij})=B_{ij}$ and
$\Psi(A)=\Psi(A_{11})+\Psi(A_{12})+\Psi(A_{22})$. Then
 we   define a map $\Psi:{\rm Alg}{\mathcal N}\rightarrow{\rm
Alg}{\mathcal M}$ and a map $h:{\rm Alg}{\mathcal N}\rightarrow
{\mathbb F}I$ with $h(A)=\Phi(A)-\Psi(A)\in{\mathbb F}I$. Now,
imitating the proof of \cite[Lemmas 2.10-2.13]{QH}, one can show
that   $\Psi: {\rm Alg}{\mathcal N}\rightarrow {\rm Alg}{\mathcal
M}$ is a ring isomorphism and  $h: {\rm Alg}{\mathcal N}\rightarrow
{\mathbb F}I$ is an additive map satisfying $h([A,B])=0$ for all
$A,B$. Hence  Claim 4.1 is true.

{\bf Claim 4.2.}  If $\Phi$ satisfies (b), then $\Phi=-\Psi+h$,
where $\Psi:{\rm Alg}{\mathcal N}\rightarrow {\rm Alg}{\mathcal M}$
is a ring anti-isomorphism and $h:{\rm Alg}{\mathcal
N}\rightarrow{\mathbb F}I$ is an additive map vanishing on all
commutators.

If (b) holds, then ${\rm ran}Q_0=(0)_+\in{\mathcal M}$ and $(0)_+$
is complemented in $\mathcal M$ by Lemmas 2.3 and 2.8. Consider the
map $\Phi': {\rm Alg}{\mathcal N}\rightarrow({\rm Alg}{\mathcal
M})^*$ defined by $\Phi'(A)=-\Phi(A)^*$ for all $A\in{\rm
Alg}{\mathcal N}$.

Imitating the proof of \cite[Lemma 2.14]{QH}, one can show that
$\Phi':{\rm Alg}{\mathcal N}\rightarrow({\rm Alg}{\mathcal M})^*$ is
a Lie multiplicative bijective map. Since, for any nontrivial
idempotent operator $P\in {\rm Alg}{\mathcal N}$,
$\Phi(P)=\tilde{\Phi}(P)+\lambda_PI$ for some $\lambda_P\in {\mathbb
C}$, we have ${\Phi}'(P)=-\tilde{\Phi}(P)^*-\lambda_PI$. Now define
a map $\tilde{{\Phi}'}:{\mathcal E}{(\mathcal
N)}\rightarrow{\mathcal E}{(\mathcal M)}^*$ by
$\tilde{{\Phi}'}(P)=I-\tilde{\Phi}(P)^*$ for all idempotents $P$.
Since $\Phi$ satisfies (b), for any nonzero idempotent $P_1\in {\rm
Alg}{\mathcal N}$, if $P_1<P_0$, we have
$\tilde{{\Phi}'}(P_1)=I-\tilde{\Phi}(P_1)^*<I-\tilde{\Phi}(P_0)^*=\tilde{{\Phi}'}(P_0)$;
if $P_1>P_0$, we have
$\tilde{{\Phi}'}(P_1)=I-\tilde{\Phi}(P_1)^*>I-\tilde{\Phi}(P_0)^*=\tilde{{\Phi}'}(P_0)$.
Hence ${\Phi}'$ satisfies (a).

Note that ${\mathcal M^\perp}=\{M^\perp : M\in{\mathcal M}\}$ is a
nest on $Y^*$. Since ${\rm ran}Q_0\in{\mathcal M}$, we have $\ker
Q_0^*=({\rm ran}Q_0)^\perp\in{\mathcal M}^\perp$, and so ${\rm
ran}(I-Q_0^*)\in{\mathcal M}^\perp $. With respect to the
decomposition $Y^*={\rm ran}(I-Q_0^*)\dot{+}{\rm ran}Q_0^*$, we have
${\rm Alg}{\mathcal M}^\perp={\mathcal D}_{11}\dot{+}{\mathcal
D}_{12}\dot{+}{\mathcal D}_{22}$, where ${\mathcal
D}_{11}=(I-Q_0^*)({\rm Alg}{\mathcal M}^\perp)(I-Q_0^*)$, ${\mathcal
D}_{12}=(I-Q_0^*)({\rm Alg}{\mathcal M}^\perp)Q_0^*$, ${\mathcal
D}_{22}=Q_0^*({\rm Alg}{\mathcal M}^\perp)Q_0^*$.  Since ${\rm
Alg}{\mathcal M}=Q_0({\rm Alg}{\mathcal M})Q_0\dot{+}Q_0({\rm
Alg}{\mathcal M})(I-Q_0)\dot{+}(I-Q_0)({\rm Alg}{\mathcal
M})(I-Q_0)={\mathcal B}_{11}\dot{+}{\mathcal B}_{12}\dot{+}{\mathcal
B}_{22}$, we have ${\mathcal B}_{11}^*\subseteq{\mathcal D}_{22}$,
${\mathcal B}_{12}^*\subseteq{\mathcal D}_{12}$ and ${\mathcal
B}_{22}^*\subseteq{\mathcal D}_{11}$. Hence $({\rm Alg}{\mathcal
M})^*={\mathcal B}_{22}^*\dot{+}{\mathcal B}_{12}^*\dot{+}{\mathcal
B}_{11}^*\subseteq{\rm Alg}{\mathcal M}^\perp$.

{\bf Step 1.} $\Phi'({\mathcal A}_{12})={\mathcal B}_{12}^*$.

The proof is the same as that of \cite[Lemma 2.15]{QH}.

{\bf Step 2.} $\Phi'({\mathcal A}_{ii})\subseteq{\mathcal
B}_{jj}^*+{\mathbb F}I$, $i,j=1,2$ and $i\not=j$.

For any $A_{11}\in{\mathcal A}_{11}$, write
$\Phi'(A_{11})=S_{22}^*+S_{12}^*+S_{11}^*$, where
$S_{ij}\in{\mathcal B}_{ij}$. Then
$$0=\Phi'([A_{11},P_0])
=[\Phi'(A_{11}),\Phi'(P_0)]=[\Phi'(A_{11}), -Q_0^*]=[Q_0^*,
\Phi'(A_{11})],$$ which implies that $S_{12}^*=0$.

Let $P\in {\mathcal A}_{22}$ be any idempotent with $P<I-P_0$. As
$\Phi'$ satisfies (a), we have
$\widetilde{\Phi'}(P)<\widetilde{\Phi'}(I-P_0)=I-
\tilde{\Phi}(P_0)^*=Q_0^*$. It follows from $[A_{11}, P]=0$ that
$0=[\Phi'(A_{11}),{\Phi}'(P)]=[S_{11}^*,\Phi'(P)]$. Since $P$ is
arbitrary,  we see that $S_{11}^*$ commutes with every idempotent in
${\mathcal B}_{11}^*$, which implies that $S_{11}$ commutes with
every idempotent in ${\mathcal B}_{11}$. Noting that
$\mathcal{B}_{11} = \mathcal{B}(\ker Q_0)$ in this case, so we must
have $S_{11}\in {\mathbb F} I_{\ker Q_0}$. By the arbitrariness of
$A_{11}$ we obtain that $\Phi'(\mathcal{A}_{11})\subseteq
\mathcal{B}_{22}^* + \mathbb{F} I$.

Similarly, one can show that, for any $A_{22}\in {\mathcal A}_{22}$,
$\Phi'(A_{22})=T_{22}^*+T_{11}^*$, where  $T_{22}$ commutes with
every idempotent in ${\mathcal B}_{22}$. Taking any $B_{22}\in
\mathcal{B}_{22}$, by the surjectivity of $\Phi'$, there exists some
$A_{0}\in {\rm Alg}\mathcal{N}$ such that $\Phi'(A_{0})=B_{22}^*$.
Furthermore, $A_{0}=A_{11}+\lambda I$ for some $A_{11}\in{\mathcal
A}_{11}$ and some scalar $\lambda$. Then
$$[B_{22}^*,T_{22}^*] =[B_{22}^*,T_{22}^*+T_{11}^*]=[\Phi'(A_{0}), \Phi'(A_{22})] = \Phi'([A_{0},A_{22}])=0$$ for all $B_{22}\in \mathcal{B}_{22}$.
Thus $[B_{22} ,T_{22} ]=0$ for all $B_{22}\in{\mathcal B}_{22}$,
which implies $T_{22}\in \mathbb{F} I_{{\rm ran}Q_0}$ as ${\mathcal
B}_{22}$ is a nest algebra. Therefore $\Phi'({\mathcal
A}_{22})\subseteq{\mathcal B}_{11}^*+{\mathbb F}I$.

By Steps 1-2, for each $A_{12}\in {\mathcal A}_{12}$, there exists
$B_{12}\in {\mathcal B}_{12}$ such that $\Phi'(A_{12})=B_{12}^*$;
for each $A_{ii}\in {\mathcal A}_{ii}$ ($i=1,2$), there exist
$B_{jj}\in {\mathcal B}_{jj}$ and $\lambda_{jj} \in {\mathbb F}$
such that $\Phi'(A_{ii})=B_{jj}^*+\lambda_{jj} I$, $i,j=1,2$ and
$i\not=j$. It is easily seen that $B_{jj}$ and $\lambda_{jj}$ are
uniquely determined. Define a map $\Psi':{\rm Alg}{\mathcal
N}\rightarrow ({\rm Alg}{\mathcal M})^*$ by
$\Psi'(A_{12})=B_{12}^*$, $\Psi'(A_{11})=B_{22}^*$,
$\Psi'(A_{22})=B_{11}^*$ and
$\Psi'(A)=\Psi'(A_{11})+\Psi'(A_{12})+\Psi'(A_{22})$ for
$A=A_{11}+A_{12}+A_{22}$.   Set $h'(A)=\Phi'(A)-\Psi'(A)$. It is
clear that $h'$ maps ${\rm Alg}{\mathcal N}$ into ${\mathbb F}I$. By
\cite[Lemma 2.17]{QH}, one can show that $\Psi'$ is a ring
isomorphism. Since, for any $A\in{\rm Alg}{\mathcal N}$, there
exists a unique element $S\in{\rm Alg}{\mathcal M}$ such that
$\Psi'(A)=S^*$, we can define a map $\Psi:{\rm Alg}{\mathcal
N}\rightarrow{\rm Alg}{\mathcal M}$ by $\Psi(A)=S$. Thus
$\Psi(A)^*=\Psi'(A)$ for every $A$ and hence $\Psi$ is a ring
anti-isomorphism.  Let $h:{\rm Alg}{\mathcal N}\rightarrow{\mathbb
F}I$ be the map defined by $h'(A)=-h(A)^*$ for all $A\in{\rm
Alg}{\mathcal N}$. Clearly, $h([A,B])=0$ for all $A,B\in{\rm
Alg}{\mathcal N}$. Furthermore, we have
$$-\Phi(A)^*=\Phi'(A)
=\Psi'(A)+h'(A)=\Psi(A)^*+h'(A) =(\Psi(A)-h(A))^*,$$ which yields
$\Phi(A)=-\Psi(A)+h(A)$ for every $A\in{\rm Alg}{\mathcal N}$. This
completes the proof of Claim 4.2.

By Claim 4.1 and Claim 4.2, Theorem 1.1 holds  for the case that
$X_-\not=X$ and $X_-$ is complemented in $X$.\hfill$\Box$

\section{ The case that $  X_{-}\not=X$ and $X_{-}$ is not complemented or $(0)\not= (0)_+ $ and $(0)_+$ is not complemented}

In this section, we deal with the case that  $  X_{-}\not=X$ and
$X_{-}$ is not complemented or $(0)\not= (0)_+ $ and $(0)_+$ is not
complemented. Here we borrow some ideas developed in \cite{WL}. Note
that, not like \cite{WL}, we do not assume that all non-trivial
elements in the nests are not complemented. Also, we give only the
detail of our proof for the
  case that $  X_{-}\not=X$ and $ X_{-}$ is not
complemented in $X$. The other case can be checked similarly.

Assume that $  X_{-}\not=X$ and $ X_{-}$ is not complemented in $X$.

Note that, there are three possible situations that $(0)$ may have,
that is, (1$^\circ$) $(0)_+\not=(0)$ and $(0)_+$ is complemented in
$X$, (2$^\circ$)  $(0)_+=(0)$ and (3$^\circ$) $(0)_+\not=(0)$ and
$(0)_+$ is not complemented in $X$. By   Section 4, Theorem 1.1 is
true if the situation (1$^\circ$) occurs. So what we need to deal
with is either (2$^\circ$) or (3$^\circ$).

Recall that $\hat{\Phi}$, $\bar{\Phi}$ and $\tilde{\Phi}$ are maps
defined in Eqs.(2.2)-(2.4), and $\Omega_1({\mathcal N},E)$ and
$\Omega_2({\mathcal N},E)$ are defined in Eq.(2.5). By Lemma 2.6,
either $\tilde{\Phi}(\Omega_1({\mathcal
N},X_-))\subseteq\Omega_1({\mathcal N},\hat{\Phi}(X_-))$ or
$\tilde{\Phi}(\Omega_1({\mathcal N},X_-))\subseteq
I-\Omega_2({\mathcal N},\hat{\Phi}(X_-))$.

The following lemma is crucial for our purpose.

{\bf Lemma 5.1.} {\it Assume that $(0)<X_-<X$ and $X_-$ is not
complemented. The following statements are true:}

(1) {\it If $\tilde{\Phi}(\Omega_1({\mathcal
N},X_-))\subseteq\Omega_1({\mathcal N},\hat{\Phi}(X_-))$, then
$\Phi(\mathbb{F}I + X\otimes X_-^{\bot})=\mathbb{F} I + Y\otimes
Y_-^{\bot}$, and there exists a ring automorphism $\tau: \mathbb{F}
\rightarrow \mathbb{F}$, a map $\gamma: X\otimes
X_-^{\bot}\rightarrow \mathbb{F}$, bijective $\tau$-linear  maps $C:
X\rightarrow Y$ and $D: X_-^{\bot} \rightarrow Y_-^{\bot}$ such that
$$\Phi(x\otimes f) =\gamma(x,f) I + Cx\otimes Df $$ holds for
all $x\in X$ and $f\in X_-^{\bot}$. }

(2) {\it  If $\tilde{\Phi}(\Omega_1({\mathcal N},X_-))\subseteq
I-\Omega_2({\mathcal N},\hat{\Phi}(X_-))$, then $\Phi(\mathbb{F}I +
X\otimes X_-^{\bot})=\mathbb{F} I+(0)_+\otimes Y^*$, and there
exists a ring automorphism $\tau: \mathbb{F} \rightarrow
\mathbb{F}$, a map $\gamma: X\otimes X_-^{\bot}\rightarrow
\mathbb{F}$, bijective $\tau$-linear  maps $C: X\rightarrow Y^*$ and
$D: X_-^{\bot} \rightarrow (0)_+$ such that
$$\Phi(x\otimes f) =\gamma(x,f) I + Df\otimes Cx  $$ holds for
all $x\in X$ and $f\in X_-^{\bot}$. }

To prove Lemma 5.1, we  consider  two cases, that is, the case that
${\mathcal N}$ has at least two nontrivial elements and the case
that ${\mathcal N}$ has only one nontrivial element. These will be
done by Lemma 5.3 and Lemma 5.4, respectively.

We first consider the case that ${\mathcal N}$ has at least two
nontrivial elements.

The following lemma is crucial for the proof of Lemma 5.3.

{\bf Lemma 5.2.} {\it Assume that $\mathcal{N}$ has at least two
nontrivial elements and $ X_{-} $ is not complemented in $X$. Then
the following statements hold.}

(i) {\it $\tilde{\Phi}(\Omega_1({\mathcal
N},X_-))\subseteq\Omega_1({\mathcal N},\hat{\Phi}(X_-))$ if and only
if $\hat{\Phi}$ is order-preserving.}

(ii) {\it $\tilde{\Phi}(\Omega_1({\mathcal N},X_-))\subseteq
I-\Omega_2({\mathcal N},\hat{\Phi}(X_-))$ if and only if
$\hat{\Phi}$ is order-reversing.}

{\bf Proof.} By Lemma 2.6, it suffices to show  that
 $\tilde{\Phi}(\Omega_1({\mathcal N},X_-))\subseteq\Omega_1({\mathcal
N},\hat{\Phi}(X_-))$ implies that $\hat{\Phi}$ is order-preserving
and $\tilde{\Phi}(\Omega_1({\mathcal N},X_-))\subseteq
I-\Omega_2({\mathcal N},\hat{\Phi}(X_-))$ implies that$\hat{\Phi}$
is order-reversing.

Assume that $\tilde{\Phi}(\Omega_1({\mathcal
N},X_-))\subseteq\Omega_1({\mathcal N},\hat{\Phi}({X_-}))$. Since
$\mathcal{N}$ has at least two nontrivial elements, we may take a
nontrivial element $E\in {\mathcal N}$ such that $E<X_-$. If, on the
contrary, $\hat{\Phi}$ is  order-reversing, then we have
$\hat{\Phi}(E)>\hat{\Phi}(X_-)$. Fix an idempotent
$P\in\Omega_{1}({\mathcal N},X_-)$. Then we have $\tilde{\Phi}(P)\in
\Omega_{1}({\mathcal M},\hat{\Phi}(X_-))$. For any $D_1\in {\mathcal
J}({\mathcal M},\hat{\Phi}(X_-))$, $D_2\in {\mathcal J}({\mathcal
M},\hat{\Phi}(E))$,  let $C_1=\bar{\Phi}^{-1}(D_1)$,
$C_2=\bar{\Phi}^{-1}(D_2)$. By the definition of $\bar{\Phi}$, we
have $C_1\in {\mathcal J}({\mathcal N},X_-)$ and $C_2\in {\mathcal
J}({\mathcal N},E)$. Furthermore,
$$\begin{array}{rl}0=&\Phi([C_{1},[C_{2},P]])=[D_{1},[D_{2},\tilde{\Phi}(P)]]\\
=&[D_{1},D_{2}\tilde{\Phi}(P)-\tilde{\Phi}(P)D_{2}]=D_{1}D_{2}\tilde{\Phi}(P)-D_{1}\tilde{\Phi}(P)D_{2}.\end{array}$$
For any $y_1\in \hat{\Phi}(X_-)$, $y_2\in\hat{\Phi}(E)$ and any
$g_1\in \hat{\Phi}(X_-)^{\bot}$, $g_2\in \hat{\Phi}(E)^{\bot}$, it
is obvious that $y_1\otimes g_1\in {\mathcal J}({\mathcal
M},\hat{\Phi}(X_-))$ and $y_2\otimes g_2\in {\mathcal J}({\mathcal
M},\hat{\Phi}(E))$. Letting $D_{i}=y_{i}\otimes g_{i}$,
 the above equation yields
$$\langle y_{2},g_{1}\rangle y_{1}\otimes
\tilde{\Phi}(P)^{*}g_{2}=\langle\tilde{\Phi}(P)y_{2},g_{1}\rangle
y_{1}\otimes g_{2}.\eqno(5.1)$$ Choose  $y_2\in\hat{\Phi}(E)$ and
$g_1\in \hat{\Phi}(X_-)^{\bot}$ such that $\langle
y_{2},g_{1}\rangle\neq0$. Eq.(5.1) implies that there exists some
scalar $\lambda_{g_2}$ such that
$\tilde{\Phi}(P)^{*}g_{2}=\lambda_{g_2} g_2$ for each $g_2\in
\hat{\Phi}(E)^{\bot}$.  It follows  that there is a scalar $\lambda$
such that $\tilde{\Phi}(P)^{*}g_{2}=\lambda g_{2}$ for all $g_{2}\in
\hat{\Phi}(E)^{\bot}$. Since $\tilde{\Phi}(P)^{*}$ is an idempotent,
either $\lambda=0$ or $\lambda=1$.

If $\lambda=0$, then $\tilde{\Phi}(P)^{*}\hat{\Phi}(E)^{\bot}=\{0\}$
and Eq.(5.1) yields $\langle\tilde{\Phi}(P)y_{2},g_{1}\rangle=0$ for
all $y_2\in\hat{\Phi}(E)$ and $g_1\in \hat{\Phi}(X_-)^{\bot}$. It
follows that $\tilde{\Phi}(P)^{*}\hat{\Phi}(X_-)^{\bot}\subseteq
\hat{\Phi}(E)^{\bot}$. So
$\tilde{\Phi}(P)^{*}\hat{\Phi}(X_-)^{\bot}=\tilde{\Phi}(P)^*(\tilde{\Phi}(P)^{*}\hat{\Phi}(X_-)^{\bot})\subseteq
\tilde{\Phi}(P)^{*}\hat{\Phi}(E)^{\bot}=\{0\}$, which is impossible
as $\tilde{\Phi}(P)\in
 \Omega_{1}({\mathcal M},\hat{\Phi}(X_-))$.

If $\lambda=1$, then Eq.(5.1) yields
$\langle\tilde{\Phi}(P)y_{2},g_{1}\rangle=\langle
y_{2},g_{1}\rangle$ for all $y_2\in\hat{\Phi}(E)$ and  $g_1\in
\hat{\Phi}(X_-)^{\bot}$. This implies that
$(I-\tilde{\Phi}(P))^{*}\hat{\Phi}(X_-)^{\bot}\subseteq
\hat{\Phi}(E)^{\bot}$, and so
$$(I-\tilde{\Phi}(P))^{*}\hat{\Phi}(X_-)^{\bot}=(I-\tilde{\Phi}(P))^{*}((I-\tilde{\Phi}(P))^{*}\hat{\Phi}(X_-)^{\bot})\subseteq
(I-\tilde{\Phi}(P))^{*}\hat{\Phi}(E)^{\bot}=\{0\}.$$ This, together
with the fact $\tilde{\Phi}(P)\hat{\Phi}(X_-)=\{0\}$ entails ${\rm
ran}\tilde{\Phi}(P)=\hat{\Phi}(X_-)$, a contradiction. Hence
$\hat{\Phi}$ is order-preserving.

Similarly, one can show that $\hat{\Phi}$ is order-reversing if
$\tilde{\Phi}(\Omega_1({\mathcal N},X_-))\subseteq
I-\Omega_2({\mathcal N},\hat{X_-})$. \hfill$\Box$

\if The following two Lemmas give a characterization of rank one
operators.

We remark that, for this case,  the authors in \cite{WL} also give a
similar characterization for  rank one operators (\cite[Lemmas
6.2-6.4]{WL}).  However, the proofs of Lemmas 6.2-6.4 depend on
Lemma 5.6 in \cite{WL}, in which the condition that $\mathcal N$
contains at least two nontrivial elements are supposed. The case
that $\mathcal N$ only contains a nontrivial element is missed.  The
following Lemmas 5.3-5.4 are generalizations of \cite[Lemmas
6.2-6.4]{WL} and make up the gap. \fi

{\bf Lemma 5.3.} {\it Assume that $\mathcal{N}$ has at least two
nontrivial elements and $(0)<X_-<X$ is not complemented in $X$. Then
for any $x\in X$ and $f\in X_-^{\bot}$, ${\Phi}(x\otimes f)$ is the
sum of a scalar and a rank one operator. Moreover, the statements
(1) and (2) of Lemma 5.1 hold.}

{\bf Proof.} We will complete the proof of the lemma by considering
three cases.

{\bf Case 1.} $x\in X_-$ and $f\in X_-^{\bot}$.

In this case, let $R=\bar{\Phi}(x\otimes f)$. Then $R\in
\mathcal{J}(\mathcal{M},\hat{\Phi}(X_-))$ and $\Phi(x\otimes
f)-\bar{\Phi}(x\otimes f)=\Phi(x\otimes f)-R\in{\mathbb F}I$. We
show that $R$ is of rank one. Assume on the contrary that ${\rm
rank}R\geq 2$. We will induce contradiction by considering two
subcases.

{\bf Subcase 1.1.}
$\tilde{\Phi}(\Omega_{1}(\mathcal{N},X_-))\subseteq\Omega_{1}(\mathcal{M},\hat{\Phi}(X_-))$.

By Lemma 5.2, $\hat{\Phi}$ is order-preserving. It is clear in this
case that we have   $(0)<Y_{-}=\hat{\Phi}(X_-)<Y$.

Since ${\rm rank}R\geq 2$,  there are two vectors $u,v\in Y\setminus
Y_-$ such that $Ru$ and $Rv$ are linearly independent. Take $h\in
Y_-^{\bot}$ such that $h(u)=1$. Then $u\otimes h\in
\Omega_{1}(\mathcal{M},\hat{\Phi}(X_-))$. Let
$B=\tilde{\Phi}^{-1}(u\otimes h)$. By Lemma 5.2(i), $B\in
\Omega_{1}({\mathcal N},X_-)$.

If $(0)_{+}=(0)\in\mathcal{N}$, then $(0)_{+}=(0)\in\mathcal{M}$ by
Lemma 2.3. Thus $(0)=\cap\{M: (0)\not=M\in{\mathcal M}\}$ and there
exists some nontrivial element $M\in\mathcal{M}$ and $g\in M^{\bot}$
such that $g(Ru)=0$ and $g(Rv)=1.$ Obviously $M<\hat{\Phi}(X_-)$,
and so $\hat{\Phi}^{-1}(M)<X_-$. Take a nonzero vector $z\in M$ and
let $A=\bar{\Phi}^{-1}(z\otimes g).$ Then $A\in
\mathcal{J}(\mathcal{N},\hat{\Phi}^{-1}(M))$, which implies
$A^{*}\hat{\Phi}^{-1}(M)^{\bot}=\{0\}$. Thus we have
$$\Phi(A(x\otimes f)B)=\Phi([A,[x\otimes f,B]])=[z\otimes
g,[R,u\otimes h]]=(z\otimes g)R(u\otimes h)=0.$$ So either
$A(x\otimes f)=0$ or $(x\otimes f)B=0.$ If $A(x\otimes f)=0,$ then
$0=\Phi([A,x\otimes f])=[z\otimes g,R]=(z\otimes g)R\not=0,$ which
is impossible; if $(x\otimes f)B=0,$ then $0=\Phi([x\otimes
f,B])=[R,u\otimes h]=R(u\otimes h)\not=0,$ which is also impossible.

If $(0)<(0)_{+}\in\mathcal{N}$, then $(0)<(0)_{+}\in\mathcal{M}$ by
Lemma 2.3. Let $M=\hat{\Phi}((0)_{+}).$ Then
$(0)<M=\hat{\Phi}((0)_{+})<\hat{\Phi}(X_-)$ as $\mathcal M$ has at
least two nontrivial elements. By Lemma 2.8, $M$ is not complemented
and thus infinite-dimensional. So, there is a vector $z\in M$ and a
functional $g\in Y^{*}$ such that $g(Ru)=0$ and $g(z)=g(Rv)=1.$ Let
$A=\tilde{\Phi}^{-1}(z\otimes g).$ As $z\otimes g\in\Omega
_2({\mathcal M},M)$, we have $A\in \Omega_{2}(\mathcal{N},(0)_{+})$
by Lemma 2.6(3). It follows that
$A^{*}\hat{\Phi}^{-1}(M)^{\bot}=\{0\}$. Now by calculating
$\Phi(A(x\otimes f)B)$ in the same way as  above, one can get a
contradiction.

So in the case that  $\hat{\Phi}$ is order-preserving, $R$ is of
rank one.

{\bf Subcase 1.2.}
$\tilde{\Phi}(\Omega_{1}(\mathcal{N},X_-))\subseteq
I-\Omega_{2}(\mathcal{M},\hat{\Phi}(X_-))$.

By Lemma 5.2(ii), in this case $\hat{\Phi}$ is order-reversing. So,
we have $(0)<(0)_{+}=\hat{\Phi}(X_-)<Y$.

If  $Y=Y_{-}$,  there is a nontrivial element $M\in{\mathcal M}$ and
vectors $u,v\in M$ such that $Ru$ and $Rv$ are linearly independent.
Let $h\in M^{\bot}$ and   $B=\bar{\Phi}^{-1}(u\otimes h).$ Then
$B\in \mathcal{J}(\mathcal{N},\hat{\Phi}^{-1}(M))$ by Lemma 2.4.
Choose $g$ in $Y^{*}$ such that $g(Ru)=0$ and $g(Rv)=1.$ Let $z\in
\hat{\Phi}(X_-)$ and $A=\tilde{\Phi}^{-1}(z\otimes g)$. By Lemma
2.6(2), $I-A\in \Omega_{1}(\mathcal{N},X_-).$  It follows from
$\hat{\Phi}(X_-)\leq M$ that $\hat{\Phi}^{-1}(M)\leq X_-$. So we get
$$\Phi(B(x\otimes f)(I-A))=\Phi([B,[x\otimes f,I-A]])=[u\otimes
h,[R,-z\otimes g]]= 0,$$ which implies that either $B(x\otimes f)=0$
or $(x\otimes f)(I-A)=0.$ If $Bx\otimes f=0$, then we get
$0=\Phi([B,x\otimes f])=[u\otimes h,R]=-Ru\otimes h\not=0$, a
contradiction. If $x\otimes f(I-A)=0$, then $0=\Phi([x\otimes
f,I-A])=[R,-z\otimes g]=z\otimes R^*g\not=0$, again a contradiction.

If $Y_{-}<Y$, there are two vectors $u,v\in Y$ such that $Ru$ and
$Rv$ are linearly independent. If $u,v\in Y_{-}$, we take $h\in
Y_-^{\bot}$ and let $B=\Phi^{-1}(u\otimes h).$ Then $B\in
\mathcal{J}(\mathcal{N},\hat{\Phi}^{-1}(Y_-))$.  Choose $g\in Y^{*}$
such that $g(Ru)=0$ and $g(Rv)=1.$ Let $z\in \hat{\Phi}(X_-)$ and
$A=\tilde{\Phi}^{-1}(z\otimes g)$. By Lemma 2.6(2), $I-A\in
\Omega_{1}(\mathcal{N},X_-).$  Still, by assumption, we have
$\hat{\Phi}(X_-)<Y_-$, and so $\hat{\Phi}^{-1}(Y_-)<X_-$. Thus by a
similar argument to that in the preceding paragraph, one can get a
contradiction. So we can
 assume that $u\not\in Y_{-}.$ Take $h\in Y_{-}^{\bot}$ such that $h(u)=1$  and let $B=\tilde{\Phi}^{-1}(u\otimes
h).$ Then $I-B\in \Omega_{2}(\mathcal{N},\hat{\Phi}^{-1}(Y_{-}))$ by
Lemma 5.2(ii). In addition, there exists some $g\in Y^{*}$ such that
$g(Ru)=0$ and $g(Rv)=1.$ Let $z=Rv$ and
$A=\tilde{\Phi}^{-1}(z\otimes g).$ Then, as $z\in
\hat{\Phi}(X_-)=(0)_{+}$, we see that $z\otimes
g\in\Omega_2({\mathcal M}, (0)_+)$ and, by Lemma 2.6 and Lemma 5.2,
$I-A\in \Omega_{1}(\mathcal{N},X_-).$ Since $\hat{\Phi}(X_-)<Y_{-}$,
we have $\hat{\Phi}^{-1}(Y_{-})<X_-$. Calculating $\Phi(B(x\otimes
f)(I-A))$, one can get a contradiction.

Hence ${\rm rank}R=1$ and $\Phi(x\otimes f)$ is the sum of a scalar
and a rank one operator, that is, the lemma is true for the case
that $x\in X_-$ and $f\in X_-^\bot$. Moreover,
$\bar{\Phi}(X_{-}\otimes X_{-}^{\bot})=\hat{\Phi}(X_-)\otimes
\hat{\Phi}(X_-)^{\bot}$.

{\bf Case 2.} $x\in X\setminus X_-$ and $f\in X_-^\perp$ with
$\langle x,f\rangle=1$.

Let  $P=x\otimes f$. Clearly, $P\in\Omega_{1}(\mathcal{N},X_{-})$.
By Case 1, we have $\bar{\Phi}(X_{-}\otimes
X_{-}^{\bot})=\hat{\Phi}(X_-)\otimes \hat{\Phi}(X_-)^{\bot}$ and
hence $ {\Phi}(X_{-}\otimes X_{-}^{\bot})\subseteq{\mathbb
F}I+\hat{\Phi}(X_-)\otimes \hat{\Phi}(X_-)^{\bot}$. For the sake of
convenience, write $\tilde{P}=\tilde{\Phi}(P)$ and
$\widehat{X_{-}}=\hat{\Phi}(X_-)$. Then
$\Phi(P)-\tilde{P}\in{\mathbb F}I$ and it suffices to show that
$\tilde{P}$ is the sum of a scalar and a rank-1 idempotent operator.

{\bf Subcase 2.1.}
$\tilde{\Phi}(\Omega_{1}(\mathcal{N},X_{-}))\subseteq\Omega_{1}(\mathcal{M},\widehat{X_{-}})$.

In this subcase $\hat{\Phi}$ is order-preserving by Lemma 5.2(i),
and hence $\widehat{X_{-}}=Y_{-}$. Note that $\Phi({\mathbb
F}I)={\mathbb F}I$ by Lemma 2.9. So, applying the fact proved in
Case 1, we have $\Phi(\mathbb{F} I+X_{-}\otimes
X_{-}^{\bot})=\mathbb{F} I+Y_{-}\otimes Y_{-}^{\bot}.$ Thus it
follows from Lemma 2.10   that there exists a ring automorphism
$\tau: \mathbb{F} \rightarrow \mathbb{F}$ and a map
$\gamma:X_{-}\times X_{-}^{\bot}\rightarrow \mathbb{F}$ such that
either
$$\Phi(y\otimes g) = \gamma(y,g)I + Cy\otimes Dg\ \ {\rm for\ all }\ \
y\in X_{-}\ \ {\rm and}\ \ g\in X_{-}^{\bot},\eqno(5.2)$$ where
 $C:X_{-} \rightarrow Y_{-}$ and
$D:X_{-}^{\bot} \rightarrow Y_{-}^{\bot}$ are two bijective
$\tau$-linear maps; or
$$\Phi(y\otimes g) = \gamma(y,g)I + Dg\otimes Cy\
\ {\rm for\ all }\ \ y\in X_{-}\ \ {\rm and}\ \ g\in
X_{-}^{\bot},\eqno(5.3)$$ where $C:X_{-} \rightarrow
 Y_{-}^{\bot}$ and $D:X_{-}^{\bot} \rightarrow Y_{-}$ are two bijective  $\tau$-linear
maps.

We first show that Eq.(5.3) can not occur. On the contrary, if
Eq.(5.3) holds,   for any $y\in X_{-}$ and $g\in X_{-}^{\bot}$, we
have $$\Phi([y\otimes g,P]) = \Phi(y\otimes P^{*}g)=
\gamma(y,P^{*}g) I + DP^{*}g\otimes Cy$$ and
$$\Phi([y\otimes g,P]) = [Dg\otimes Cy,\tilde{P}] = Dg\otimes
\tilde{P}^{*}Cy,$$ which imply  $Dg\otimes
\tilde{P}^{*}Cy=DP^{*}g\otimes Cy$  for all $y\in X_{-}$ and $g\in
X_{-}^{\bot}$. Thus there exists some scalar $\lambda$ such that
$D|_{X_-^\bot}=\lambda DP^*|_{X_-^\bot}$. Since $P$ is of rank one,
we see that $D|_{X_-^\bot}$ is also of rank one, but this is
impossible as $X_{-}^{\bot}$ is infinite-dimensional.

So Eq.(5.2) holds. Then, for any $y\in X_{-}$ and $g\in
X_{-}^{\bot}$, we have
$$\Phi([y\otimes g,P]) = \Phi(y\otimes P^{*}g)= \gamma(y,P^{*}g) I +
Cy\otimes DP^{*}g$$ and
$$\Phi([y\otimes g,P]) = [Cy\otimes Dg,\tilde{P}] = Cy\otimes
\tilde{P}^{*}Dg.$$ It follows that $Cy\otimes \tilde{P}^{*}Dg=
\gamma(y,P^{*}g) I + Cy\otimes DP^{*}g$ holds for all $y\in X_{-}$
and $g\in X_{-}^{\bot}.$ Since $I$ is of infinite rank, we have
$\gamma(y,P^{*}g)=0$, and $Cy\otimes \tilde{P}^{*}Dg= Cy\otimes
DP^{*}g$. So $\tilde{P}^{*}Dg=DP^{*}g$ for all $g\in X_{-}^{\bot}.$
Since $P$ is of rank one, it follows that the restriction of
$\tilde{P}^{*}$ to $Y_{-}^{\bot}$ is of rank one. Note that
$\tilde{P}^{*}Y^{*} \subseteq Y_{-}^{\bot}$ and $\tilde{P}^{*}Y^{*}
= \tilde{P}^{*}(\tilde{P}^{*}Y^{*})\subseteq
\tilde{P}^{*}Y_{-}^{\bot}.$ So $\tilde{P}^{*}$ is of rank one, which
implies that  $\tilde{P}$ is also of rank one, as desired.

{\bf Subcase 2.2.}  $\tilde{\Phi}(\Omega_{1}(\mathcal{N},X_{-}))
\subseteq I-\Omega_{2}(\mathcal{M},\widehat{X_{-}})$.

By Lemma 5.2(ii),  $\hat{\Phi}$ is order-reversing, and
$\widehat{X_{-}}=(0)_{+}$.  Applying the fact proved in Case 1, we
have $\Phi(\mathbb{F} I+X_{-}\otimes X_{-}^{\bot})=\mathbb{F}
I+(0)_{+}\otimes (0)_{+}^{\bot}.$ Thus, by Lemma 2.10 again, there
exists a ring automorphism $\tau: \mathbb{F} \rightarrow \mathbb{F}$
and a map $\gamma:X_{-}\times X_{-}^{\bot}\rightarrow \mathbb{F}$
such that either
$$\Phi(x\otimes f)=\gamma(x,f)I + Cx\otimes Df\ \ {\rm for\ all }\ \
x\in X_{-}\ \ {\rm and}\ \ f\in X_{-}^{\bot},\eqno(5.4)$$ where
 $C:X_{-} \rightarrow (0)_+$ and
$D:X_{-}^{\bot} \rightarrow (0)_+^{\bot}$ are two $\tau$-linear
bijective maps; or
$$\Phi(x\otimes f)=\gamma(x,f)I + Df\otimes Cx\
\ {\rm for\ all }\ \ x\in X_{-}\ \ {\rm and}\ \ f\in
X_{-}^{\bot},\eqno(5.5)$$ where   $C:X_{-} \rightarrow
 (0)_+^{\bot}$ and $D:X_{-}^{\bot} \rightarrow (0)_+$ are two $\tau$-linear
bijective maps.

If Eq.(5.4) holds, by calculating $\Phi([y\otimes g,P])$, one
obtains $(I-\tilde{P})Cy\otimes Dg= Cy\otimes DP^{*}g$ for all $y\in
X_{-}$ and $g\in X_{-}^{\bot}.$ Since $P$ is of rank one,  we get
$D$ is also of rank one, which is impossible. So we must have that
Eq.(5.5) holds.  By calculating $\Phi([y\otimes g,P])$, one can get
$(I-\tilde{P})Dg= DP^{*}g$ for all $f\in X_{-}^{\bot}$, which
implies that the restriction of $I-\tilde{P}$ to $(0)_{+}$ is of
rank one. So $I-\tilde{P}$ is of rank one and $\tilde{P}$ is the sum
of a scalar and a rank one operator, as desired.

Summing up, we have proved that $\Phi(x\otimes f)=\Phi(P)$ is the
sum of a scalar and a rank one operator if $x\in X, f\in X_-^\bot$
with $\langle x,f\rangle=1$. Moreover,
$\tilde{\Phi}(\Omega_{1}(\mathcal{N},X_{-}))\subseteq\Omega_{1}(\mathcal{M},\widehat{X_{-}})$
implies that Eq.(5.2) holds, while
$\tilde{\Phi}(\Omega_{1}(\mathcal{N},X_{-})) \subseteq
I-\Omega_{2}(\mathcal{M},\widehat{X_{-}})$ implies that Eq.(5.5)
holds.

{\bf Case 3.} $x\in X\setminus X_-$ and $f\in X_{-}^{\bot}.$

Now assume that $x\in X\setminus X_-$ and $f\in X_{-}^{\bot}.$ We
need still consider two subcases.

{\bf Subcase 3.1.} $\langle x,f\rangle=\lambda \neq 0$.

Let $P=\lambda^{-1}x\otimes f$.  Then the  rank-one idempotent $P\in
\Omega_{1}(\mathcal{N},X_{-})$ and $x\otimes f=\lambda P.$ Write
$\tilde{P}=\tilde{\Phi}(P).$ By what proved in Case 2, $\tilde{P}$
is a rank-1 idempotent.

{\bf Subcase 3.1.1.}
$\tilde{\Phi}(\Omega_{1}(\mathcal{N},X_{-}))\subseteq\Omega_{1}(\mathcal{M},\widehat{X_{-}})$.

By Lemma 5.2(i), $\hat{\Phi}$ is order-preserving, and then
$\widehat{X_{-}}=Y_{-}$. So we still have that either Eq.(5.2) or
Eq.(5.3) holds.

If Eq.(5.2) holds, then for any $y\in X_{-}$ and $g\in
X_{-}^{\bot}$, we have
$$\Phi([y\otimes g,\lambda P]) = [\Phi(y\otimes g),\Phi(\lambda P)]=Cy\otimes \Phi(\lambda P)^{*}Dg
-\Phi(\lambda P)Cy\otimes Dg$$ and
$$\Phi([y\otimes g,\lambda P]) = [\tau(\lambda)Cy\otimes Dg,\tilde{P}] = \tau(\lambda)Cy\otimes
\tilde{P}^{*}Dg.$$ It follows that $\tau(\lambda)Cy\otimes
\tilde{P}^{*}Dg=Cy\otimes \Phi(\lambda P)^{*}Dg -\Phi(\lambda
P)Cy\otimes Dg$, that is,
$$Cy\otimes
(\overline{\tau(\lambda)}\tilde{P}^{*}Dg- \Phi(\lambda P)^{*}Dg) =
-\Phi(\lambda P)Cy\otimes Dg\ \ {\rm for \ all}\ y\in X_-\ {\rm
and}\ g\in X_-^\perp.$$ Since $C$ is bijective, there exists some
scalar $\alpha$ such that $\tau(\lambda)\tilde{P}^{*}Dg-
\Phi(\lambda P)^{*}Dg= \alpha Dg$ for all $g\in X_{-}^{\bot}$, that
is, $(\Phi(\lambda P)^{*} + \alpha I)|_{Y_{-}^{\bot}} =
\tau(\lambda)\tilde{P}^{*}|_{Y_{-}^{\bot}}$ as
$D:X_-^\bot\rightarrow Y_-^\bot$ is bijective. Let
$\widetilde{\Phi(\lambda P)} = \Phi(\lambda P) + \overline{\alpha
}I.$ Note that $[\lambda P,A] = [\lambda P,[P,[P,A]]]$ holds for all
$A\in {\rm Alg}\mathcal{N}$. So we get
$$[\widetilde{\Phi(\lambda P)},\Phi(A)] = [\widetilde{\Phi(\lambda
P)},[\tilde{P},[\tilde{P},\Phi(A)]]], \ \ \ \forall A\in {\rm
Alg}\mathcal{N}.\eqno(5.6)$$  By the bijectivity of $\Phi$, for any
$y\in Y$ and $ g\in Y_{-}^{\bot}$, there exists some $A\in{\rm
Alg}{\mathcal N}$ such that $\Phi(A) =y\otimes g$. Thus Eq.(5.6)
entails
$$(\widetilde{\Phi(\lambda P)}-\widetilde{\Phi(\lambda
P)}\tilde{P})y\otimes g= (\tau(\lambda)\tilde{P} +
\widetilde{\Phi(\lambda P)}-2\widetilde{\Phi(\lambda
P)}\tilde{P})y\otimes \tilde{P}^{*}g \eqno(5.7)$$ for all $y\in Y$
and $ g\in Y_{-}^{\bot}$. By Case 2, we can write $\tilde{P} =
u\otimes h$, where $u\in Y$ and $h\in Y_{-}^{\bot}$ with $\langle u
,h\rangle = 1.$ Since $\dim(Y_{-}^{\bot}) > 2,$ there exists some
$g_{1}\in Y_{-}^{\bot}$ such that $g_1$ is linearly independent of
$h$. If $\langle u ,g_1\rangle \neq0,$ let $g=g_{1}.$ If $\langle u
,g_1\rangle=0$, let $g=h+ g_{1}$. Then $g$ and $h$ are linearly
independent and $\langle u,g\rangle \neq0.$ So $g$ and
$\tilde{P}^{*}g$  are also linearly independent. By Eq.(5.7), we get
$(\widetilde{\Phi(\lambda P)}- \widetilde{\Phi(\lambda
P)}\tilde{P})y= 0$ for all $y\in Y.$ It follows that
$\widetilde{\Phi(\lambda P)} = \widetilde{\Phi(\lambda
P)}\tilde{P}$, which implies $\widetilde{\Phi(\lambda P)}$ is of
rank one. So $\Phi(\lambda P)=\widetilde{\Phi(\lambda P)}-\bar\alpha
I$ is the sum of a scalar and a rank one operator.

We claim that Eq.(5.3) can not occur. If, on the contrary,  Eq.(5.3)
holds, then for any $y\in X_{-}$ and $g\in X_{-}^{\bot}$, we have
$$\Phi([y\otimes g,\lambda P])=[\Phi(y\otimes g),\Phi(\lambda P)]=Dg\otimes \Phi(\lambda
P)^{*}Cy-\Phi(\lambda P)Dg\otimes Cy$$ and
$$\Phi([y\otimes g,\lambda P])=\Phi(\lambda y\otimes P^*g)=\tau(\lambda)D(P^*g)\otimes Cy.$$ It follows that
$Dg\otimes\Phi(\lambda P)^{*}Cy=(\Phi(\lambda
P)Dg+\tau(\lambda)D(P^*g))\otimes Cy$ for  all $y\in X_-$ and $g\in
X_-^\perp.$ So there exists some scalar $\gamma$ such that
$\Phi(\lambda P)^{*}Cy=\gamma Cy$, which implies $\Phi(\lambda
P)^{*}=\gamma I$ on $Y_-^\perp$ as $C$ is bijective. Now, for any
$A\in {\rm Alg}\mathcal{N}$, by the relation $[\lambda P,A] =
[\lambda P,[P,[P,A]]]$, we get $[{\Phi(\lambda P)},\Phi(A)] =
[{\Phi(\lambda P)},[\tilde{P},[\tilde{P},\Phi(A)]]]$. Particularly,
for any $y\in Y$ and $g\in Y_{-}^{\bot}$,  there exists some
$A\in{\rm Alg}{\mathcal N}$ such that $\Phi(A) = y\otimes g$. Thus
we have $[{\Phi(\lambda P)},y\otimes g] = [{\Phi(\lambda
P)},[\tilde{P},[\tilde{P},y\otimes g]]]$, that is,
$$({\Phi(\lambda P)}-{\Phi(\lambda
P)}\tilde{P}-\gamma I+\gamma\tilde{P})y\otimes g= ({\Phi(\lambda
P)}-2{\Phi(\lambda P)}\tilde{P}+2\gamma\tilde{P}-\gamma I)y\otimes
\tilde{P}^{*}g $$ holds for all $y\in Y$ and $g\in Y_{-}^{\bot}$.
Still, we can choose $g$ such  that $g$ and $\tilde{P}^{*}g$  are
linearly independent. The above equation yields $({\Phi(\lambda
P)}-{\Phi(\lambda P)}\tilde{P}-\gamma
I+\gamma\tilde{P})y=({\Phi(\lambda P)}-2{\Phi(\lambda
P)}\tilde{P}+2\gamma\tilde{P}-\gamma I)y=0$ for all $y\in Y$. This
implies   ${\Phi(\lambda P)}=\gamma I$, which is contradicting to
$\Phi({\mathbb F}I)={\mathbb F}I$ and the bijectivity of $\Phi$. So
Eq.(5.3) can not occur.

{\bf Subcase 3.1.2.} $\tilde{\Phi}(\Omega_{1}(\mathcal{N},X_{-}))
\subseteq I-\Omega_{2}(\mathcal{M},\widehat{X_{-}})$

By Lemma 5.2(ii),   $\hat{\Phi}$ is order-reversing and thus
$\widehat{X_{-}}=(0)_{+}$.   It follows that either Eq.(5.4) or
Eq.(5.5) holds. By a similar argument to the Subcase 3.1.1 above,
one can check that, Eq.(5.4) can not occur and if Eq.(5.5) holds,
then $\Phi(x\otimes f)=\Phi(\lambda P)$ is the sum of a scalar and a
rank one operator.

{\bf Subcase 3.2.} $\langle x,f\rangle=0$.

In this case,  take $x_{1}\in X$ such that $\langle x_1,f\rangle=1$
and let $x_{2}=x-2x_{1}$, $x_{3}=x - x_{1}.$ Then $\langle
x_{i},f\rangle \neq0$ for $i=1,2,3$. By Subcase 3.1,
$\Phi(x_i\otimes f)$ is the sum of a scalar and a rank one operator.
So we can assume $\Phi(x_{i}\otimes f) = \lambda_{i} I +
u_{i}\otimes h_{i}$ for some $\lambda_i$, $i=1,2,3.$ Note that
$$(\lambda_{1} + \lambda_{2})I + u_{1}\otimes h_{1} + u_{2}\otimes
h_{2} = \Phi((x_{1} + x_{2})\otimes f) = \Phi(x_{3}\otimes f) =
\lambda_{3} I + u_{3}\otimes h_{3}.$$ Since $I$ is of infinite rank,
we must have $u_{1}\otimes h_{1} + u_{2}\otimes h_{2} = u_{3}\otimes
h_{3}$. This forces that $\{u_{1},u_{3}\}$ or $\{h_{1},h_{3}\}$ is a
linearly dependent set. So $u_{1}\otimes h_{1} + u_{3}\otimes h_{3}$
is of rank one. Thus $\Phi(x\otimes f) = \Phi( (x_{1} +
x_{3})\otimes f ) = (\lambda_{1} + \lambda_{3})I + u_{1}\otimes
h_{1} + u_{3}\otimes h_{3}$ is the sum of a scalar and a rank one
operator.

Combining Cases 1-3 and the bijectivity of $\Phi$, we have shown
that

(i) $\Phi({\mathbb F}I+X\otimes X_-^\bot)={\mathbb F}I+Y\otimes
Y_-^\bot$ if
$\tilde{\Phi}(\Omega_{1}(\mathcal{N},X_{-}))\subseteq\Omega_{1}(\mathcal{M},\widehat{X_{-}})$;

(ii) $\Phi({\mathbb F}I+X\otimes X_-^\bot)={\mathbb F}I+(0)_+\otimes
Y^*$ if $\tilde{\Phi}(\Omega_{1}(\mathcal{N},X_{-})) \subseteq
I-\Omega_{2}(\mathcal{M},\widehat{X_{-}})$.

So Lemma 2.10 is applicable. Observing the arguments in Case 2 and
Case 3, it is easily seen that (i) implies Lemma 5.1(1) and (ii)
implies Lemma 5.1(2).

The proof of  Lemma 5.3 is finished. \hfill$\Box$

For the case that the nest has only one nontrivial element, we have

{\bf Lemma 5.4.} {\it Assume that $\mathcal{N}=\{(0),X_-,X\}$ and
$X_-$ is not complemented in $X$. Then for any $x\in X$ and $f\in
X_-^{\bot}$, ${\Phi}(x\otimes f)$ is the sum of a scalar and a rank
one operator. Moreover, the statements (1) and (2) of Lemma 5.1
hold.}

{\bf Proof.} Not that $(0)_+=X_-$ in this situation. By Lemma 2.3,
${\mathcal M}=\{(0),Y_-,Y\}$ and $(0)_+=Y_-$ is not complemented in
$Y$. Let  $R\in \mathcal{J}(\mathcal{M},Y_-)$ so that
${\Phi}(x\otimes f)-R\in{\mathbb F}I$ (Lemma 2.2). We remak here
that, since $\mathcal N$ only contains one nontrivial element, Lemma
5.2 and Lemma 2.6(3) are not  applicable. As in the proof of Lemma
5.3, we consider three cases.

{\bf Case 1.} $x\in X_-$ and $f\in X_-^{\bot}$.

In this case, $x\otimes f\in{\mathcal J}({\mathcal N}, X_-)$.
We'll prove ${\rm
rank}R=1$. Assume on the contrary that ${\rm rank}R\geq 2$. Then
there are two vectors $u,v\in Y\setminus Y_-$ such that $Ru$ and
$Rv$ are linearly independent. Take $h\in Y_-^{\bot}$ such that
$h(u)=1$. Then $u\otimes h\in \Omega_{1}(\mathcal{M},Y_-)$. Note
that $Y_-$ is infinite-dimensional. There is a vector $z\in Y_-$ and
a functional $g\in Y^{*}$ such that $g(Ru)=0$ and $g(z)=g(Rv)=1$.
Let $A=\tilde{\Phi}^{-1}(z\otimes g)$ and
$B=\tilde{\Phi}^{-1}(u\otimes h)$.

By Lemma 2.6(1), either
$\tilde{\Phi}(\Omega_{1}(\mathcal{N},X_-))\subseteq\Omega_{1}(\mathcal{M},Y_-)$
or
$I-\tilde{\Phi}(\Omega_{1}(\mathcal{N},X_-))\subseteq\Omega_2(\mathcal{M},Y_-)$.
If
$\tilde{\Phi}(\Omega_{1}(\mathcal{N},X_-))\subseteq\Omega_{1}(\mathcal{M},Y_-)$,
then $B\in \Omega_{1}({\mathcal N},X_-)$, and by Lemma 2.6(2), $A\in
\Omega_{2}(\mathcal{N},X_-)$. Thus we get $\Phi(A(x\otimes
f)B)=\Phi([A,[x\otimes f,B]])=[z\otimes g,[R,u\otimes h]]=(z\otimes
g)R(u\otimes h)=0$, which implies either $A(x\otimes f)=0$ or
$(x\otimes f)B=0.$ If $A(x\otimes f)=0,$ then $0=\Phi([A,x\otimes
f])=[z\otimes g,R]=(z\otimes g)R,$ which is impossible; if
$(x\otimes f)B=0,$ then $0=\Phi([x\otimes f,B])=[R,u\otimes
h]=R(u\otimes h),$ which is also impossible. If
$I-\tilde{\Phi}(\Omega_{1}(\mathcal{N},X_-))\subseteq\Omega_2(\mathcal{M},Y_-)$,
then $I-B\in \Omega_{2}(\mathcal{N},X_-)$, and by Lemma 2.6(2),
$I-A\in \Omega_{1}(\mathcal{N},X_-).$  Thus $\Phi((I-B)(x\otimes
f)(I-A))=\Phi([I-B,[x\otimes f,I-A]])=[-u\otimes h,[R,-z\otimes
g]]=0$. So we get either $(I-B)(x\otimes f)=0$ or $(x\otimes
f)(I-A)=0.$ Still, this is impossible.

Hence ${\rm rank}R=1$.

{\bf Case 2.} $x\in X\setminus X_-$ and $f\in X_-^\perp$ with
$\langle x,f\rangle=1$.

In this case $P=x\otimes f\in\Omega_{1}(\mathcal{N},X_{-})$. Write
$\tilde{P}=\tilde{\Phi}(P)$.  By Case 1 and Lemma 2.9, we have
  $\Phi(\mathbb{F}
I+X_{-}\otimes X_{-}^{\bot})=\mathbb{F} I+Y_{-}\otimes
Y_{-}^{\bot}.$ So, by Lemma 2.10, there exists a ring automorphism
$\tau: \mathbb{F} \rightarrow \mathbb{F}$ and a map
$\gamma:X_{-}\times X_{-}^{\bot}\rightarrow \mathbb{F}$ such that
either
$$\Phi(y\otimes g) = \gamma(y,g)I + Cy\otimes Dg\ \ {\rm for\ all }\ \
y\in X_{-}\ \ {\rm and}\ \ g\in X_{-}^{\bot},\eqno(5.8)$$ where
 $C:X_{-} \rightarrow Y_{-}$ and
$D:X_{-}^{\bot} \rightarrow Y_{-}^{\bot}$ are two $\tau$-linear
bijective maps; or
$$\Phi(y\otimes g) = \gamma(y,g)I + Dg\otimes Cy\
\ {\rm for\ all }\ \ y\in X_{-}\ \ {\rm and}\ \ g\in
X_{-}^{\bot},\eqno(5.9)$$ where $C:X_{-} \rightarrow
 Y_{-}^{\bot}$ and $D:X_{-}^{\bot} \rightarrow Y_{-}$ are two $\tau$-linear
bijective maps.

Assume first that  Eq.(5.8) holds. If
$\tilde{\Phi}(\Omega_{1}(\mathcal{N},X_-))\subseteq\Omega_{1}(\mathcal{M},Y_-)$,
then $\tilde{P}\in \Omega_{1}({\mathcal M},Y_-)$. Thus, for any
$y\in Y$ and any $g\in Y_-^\perp$, we have $$\Phi([y\otimes g,P]) =
\Phi(y\otimes P^{*}g)= \gamma(y,P^{*}g) I + Cy\otimes DP^{*}g$$ and
$$\Phi([y\otimes g,P]) = [Cy\otimes Dg,\tilde{P}] = Cy\otimes
\tilde{P}^{*}Dg.$$ Since $I$ is of infinite rank, the above two
equations yield  $Cy\otimes \tilde{P}^{*}Dg=Cy\otimes DP^{*}g$, and
so $\tilde{P}^{*}Dg =DP^{*}g$ for all $g\in X_{-}^{\bot}.$ Since $P$
is of rank one,  the restriction of $\tilde{P}^{*}$ to
$Y_{-}^{\bot}$ is of rank one. Note that $\tilde{P}^{*}Y^{*}
\subseteq Y_{-}^{\bot}$ and $\tilde{P}^{*}Y^{*} =
\tilde{P}^{*}(\tilde{P}^{*}Y^{*})\subseteq
\tilde{P}^{*}Y_{-}^{\bot}.$ So $\tilde{P}^{*}$ is of rank one, which
implies that  $\tilde{P}$ is also of rank one. Hence $\Phi(P)$ is
the sum of a scalar and a rank one operator.

If
$I-\tilde{\Phi}(\Omega_{1}(\mathcal{N},X_-))\subseteq\Omega_2(\mathcal{M},Y_-)$,
then $I-\tilde{P}\in \Omega_{2}(\mathcal{M},Y_-)$. For any $y\in Y$
and any $g\in Y_-^\perp$, we have $$\Phi([y\otimes g,P]) =
\Phi(y\otimes P^{*}g)= \gamma(y,P^{*}g) I + Cy\otimes DP^{*}g$$ and
$$\Phi([y\otimes g,P])=[Cy\otimes Dg,\tilde{P}]=-[Cy\otimes Dg,I-\tilde{P}]=(I-\tilde{P})Cy\otimes
Dg.$$ The above two equations yield $(I-\tilde{P})Cy\otimes
Dg=Cy\otimes DP^{*}g$ for all $y\in X_-$ and $g\in X_{-}^{\bot}$,
and hence $D$ and $DP^*$ are linearly dependent. Since $P$ is of
rank one, $D$ is also of rank one, which is impossible as
$X_{-}^{\bot}$ is infinite-dimensional. Therefore, this case can not
occur.

Similarly, if Eq.(5.9) holds, one can show that
$I-\tilde{\Phi}(\Omega_{1}(\mathcal{N},X_-))\subseteq\Omega_2(\mathcal{M},Y_-)$
and $\Phi(P)$ is the sum of a scalar and a rank one operator.

{\bf Case 3.} $x\in X\setminus X_-$ and $f\in X_{-}^{\bot}.$

Note that, we still have that either  Eq.(5.8) or Eq.(5.9) holds.

{\bf Subcase 3.1.} $\langle x,f\rangle=\lambda \neq 0$.

Then there exists a rank-one idempotent $P\in
\Omega_{1}(\mathcal{N},X_{-})$ such that $x\otimes f=\lambda P.$
Write $\tilde{P}=\tilde{\Phi}(P)$. By Case 2, $\tilde{P}$ is a
rank-1 idempotent.

Assume that  Eq.(5.8) holds. If
$\tilde{\Phi}(\Omega_{1}(\mathcal{N},X_-))\subseteq\Omega_{1}(\mathcal{M},Y_-)$,
then $\tilde{P}\in \Omega_{1}({\mathcal M},Y_-)$.  So, for any $y\in
X_{-}$ and $g\in X_{-}^{\bot}$, we have  $$\Phi([y\otimes g,\lambda
P])=[\Phi(y\otimes g),\Phi(\lambda P)]=Cy\otimes \Phi(\lambda
P)^{*}Dg -\Phi(\lambda P)Cy\otimes Dg$$ and
$$\Phi([y\otimes g,\lambda P])=\Phi([\lambda y\otimes g,P])=[\tau(\lambda)Cy\otimes Dg,\tilde{P}] = \tau(\lambda)Cy\otimes
\tilde{P}^{*}Dg.$$ It follows that
$$Cy\otimes
(\overline{\tau(\lambda)}\tilde{P}^{*}Dg- \Phi(\lambda P)^{*}Dg) =
-\Phi(\lambda P)Cy\otimes Dg\ {\rm for \ all}\ y\in X_-\ {\rm and}\
g\in X_-^\perp.$$ Hence there exists some scalar $\alpha$ such that
$\tau(\lambda)\tilde{P}^{*}Dg- \Phi(\lambda P)^{*}Dg= \alpha Dg$ for
all $g\in X_{-}^{\bot}$, which implies that $(\Phi(\lambda P)^{*} +
\alpha I)|_{Y_{-}^{\bot}}=
\tau(\lambda)\tilde{P}^{*}|_{Y_{-}^{\bot}}$ as $D$ is bijective. Let
$\widetilde{\Phi(\lambda P)} = \Phi(\lambda P) + \overline{\alpha}
I.$ For any $A\in {\rm Alg}\mathcal{N}$, by the relation  $[\lambda
P,A] = [\lambda P,[P,[P,A]]]$, we get $[\widetilde{\Phi(\lambda
P)},\Phi(A)] = [\widetilde{\Phi(\lambda
P)},[\tilde{P},[\tilde{P},\Phi(A)]]]$. Particularly, for any $y\in
Y$ and $g\in Y_{-}^{\bot}$,  by the bijectivity of $\Phi$,  there
exists some $A\in{\rm Alg}{\mathcal N}$ such that $\Phi(A) =
y\otimes g$. So we have $[\widetilde{\Phi(\lambda P)},y\otimes g] =
[\widetilde{\Phi(\lambda P)},[\tilde{P},[\tilde{P},y\otimes g]]]$,
that is,
$$(\widetilde{\Phi(\lambda P)} - \widetilde{\Phi(\lambda
P)}\tilde{P})y\otimes g= (\tau(\lambda)\tilde{P} +
\widetilde{\Phi(\lambda P)}-2\widetilde{\Phi(\lambda
P)}\tilde{P})y\otimes \tilde{P}^{*}g\eqno(5.10)$$ holds for all
$y\in Y$ and $g\in Y_{-}^{\bot}$. By Case 2, we can write $\tilde{P}
= u\otimes h$, where $u\in Y$ and $h\in Y_{-}^{\bot}$ with $\langle
u ,h\rangle = 1.$ Since $\dim(Y_{-}^{\bot}) > 2,$ there exists $g$
with $\langle u,g\rangle \neq0$ such that $g$ and $h$ are linearly
independent. So $g$ and $\tilde{P}^{*}g$  are also linearly
independent. By Eq.(5.10), we get $(\widetilde{\Phi(\lambda P)}-
\widetilde{\Phi(\lambda P)}\tilde{P})y= 0$ for all $y\in Y$. So
$\widetilde{\Phi(\lambda P)} = \widetilde{\Phi(\lambda
P)}\tilde{P}$, which implies $\widetilde{\Phi(\lambda P)}$ is of
rank one. So $\Phi(\lambda P)=\widetilde{\Phi(\lambda P)}-\alpha I$
is the sum of a scalar and a rank one operator.

We claim that the case
$I-\tilde{\Phi}(\Omega_{1}(\mathcal{N},X_-))\subseteq\Omega_2(\mathcal{M},Y_-)$
does not happen. Otherwise, we have  $I-\tilde{P}\in
\Omega_{2}(\mathcal{M},Y_-)$. Then for any $y\in X_{-}$ and $g\in
X_{-}^{\bot}$, by calculating $\Phi([y\otimes g,\lambda P])$, one
can obtain
$$Cy\otimes\Phi(\lambda P)^{*}Dg =
(\Phi(\lambda P)-\tau(\lambda)(I-\tilde{P}))Cy\otimes Dg\ {\rm for \
all}\ y\in X_-\ {\rm and}\ g\in X_-^\perp.$$ It follows that there
exists some scalar $\gamma$ such that $\Phi(\lambda P)^*Dg=\gamma
Dg$, which implies $\Phi(\lambda P)^*=\gamma I$ on $Y_-^\perp$ as
$D$ is bijective. Now, for any $A\in {\rm Alg}\mathcal{N}$, by the
relation  $[\lambda P,A] = [\lambda P,[P,[P,A]]]$, we get
$[{\Phi(\lambda P)},\Phi(A)] = [{\Phi(\lambda
P)},[\tilde{P},[\tilde{P},\Phi(A)]]]$. Particularly, for any $y\in
Y$ and $g\in Y_{-}^{\bot}$,  there exists some $A\in{\rm
Alg}{\mathcal N}$ such that $\Phi(A) = y\otimes g$. So we have
$[{\Phi(\lambda P)},y\otimes g] = [{\Phi(\lambda
P)},[\tilde{P},[\tilde{P},y\otimes g]]]$, that is,
$$({\Phi(\lambda P)}-{\Phi(\lambda
P)}\tilde{P}-\gamma I+\gamma\tilde{P})y\otimes g= ({\Phi(\lambda
P)}-2{\Phi(\lambda P)}\tilde{P}+2\gamma\tilde{P}-\gamma I)y\otimes
\tilde{P}^{*}g\eqno(5.11)$$ holds for all $y\in Y$ and $g\in
Y_{-}^{\bot}$. Still, we can choose $g$ such  that $g$ and
$\tilde{P}^{*}g$  are  linearly independent. By Eq.(5.11), we get
$({\Phi(\lambda P)}-{\Phi(\lambda P)}\tilde{P}-\gamma
I+\gamma\tilde{P})y=({\Phi(\lambda P)}-2{\Phi(\lambda
P)}\tilde{P}+2\gamma\tilde{P}-\gamma I)y=0$ for all $y\in Y$. This
leads to a contradiction   ${\Phi(\lambda P)}=\gamma I$.

If Eq.(5.9) holds,  by a similar argument to  the above, one can
show that
$I-\tilde{\Phi}(\Omega_{1}(\mathcal{N},X_-))\subseteq\Omega_2(\mathcal{M},Y_-)$
and $\Phi(\lambda P)$ is also the sum of a scalar and a rank one
operator.

{\bf Subcase 3.2.} $\langle x,f\rangle=0$.

The proof is the same to that of Subcase 3.2 in Lemma 5.3. We omit
it here.

Combining Cases 1-3, we see that the statements (i) and (ii) in the
proof of Lemma 5.3 still hold, and this   completes the proof of
Lemma 5.4. \hfill$\Box$

{\bf Proof of Lemma 5.1.} It is immediate from Lemma 5.3 and Lemma
5.4. \hfill$\Box$ \vspace{3mm}

 {\bf Proof of Theorem 1.1 for the
case $X_-\not=X$ and $X_-$ is not complemented}.

Now let us show that Theorem 1.1 is true  for the case that $
X_{-}\not=X$ and $X_{-}$ is not complemented. By Lemma 5.1, (1) or
(2) holds.

If Lemma 5.1(1) holds, then $\Phi(\mathbb{F}I + X\otimes
X_-^{\bot})=\mathbb{F} I + Y\otimes Y_-^{\bot}$, and there exists a
ring automorphism $\tau: \mathbb{F} \rightarrow \mathbb{F}$ and a
map $\gamma: X\otimes X_-^{\bot}\rightarrow \mathbb{F}$ such that
$$\Phi(x\otimes f) =\gamma(x,f) I + Cx\otimes Df \eqno (5.12)$$ holds for
all $x\in X$ and $f\in X_-^{\bot}$, where   $C: X\rightarrow Y$ and
$D: X_-^{\bot} \rightarrow Y_-^{\bot}$ are two $\tau$-linear
bijective maps. Thus for any $A\in{\rm Alg}{\mathcal N}$, any $x\in
X$ and any $f\in X_-^\perp$, as $\Phi([A,x\otimes
f])=[\Phi(A),\Phi(x\otimes f)]$  and  $I$ is of infinite rank, one
  obtains that
$$Cx\otimes (\Phi(A)^{*}Df- DA^{*}f )= (\Phi(A)Cx -CAx)\otimes Df.$$ As $D$ is injective, the above equation entails that there exists a scalar $h(A)$ such that
$$\Phi(A)C=CA+h(A)C. \eqno(5.13)$$ It is clear that $h$ is an additive functional
on Alg${\mathcal N}$. Define $\Psi: {\rm  Alg}{\mathcal
N}\rightarrow{\rm Alg}{\mathcal M}$ by $\Psi(A)=\Phi(A)-h(A)I$ for
all $A\in{\rm Alg}{\mathcal N}$. Then, $\Psi$ is an additive
bijection. By Eq.(5.13), for any $A,B\in{\rm Alg}{\mathcal N}$, we
have
$$\Psi(AB)C=CAB=\Psi(A)CB=\Psi(A)\Psi(B)C.$$
Since ran$C=Y$, we see that $\Psi(AB)= \Psi(A)\Psi(B)$ for all
$A,B\in {\rm Alg}\mathcal{N}$, that is, $\Psi$ is a ring
isomorphism. 

Assume that Lemma 5.1(2) holds. Then $\Phi(\mathbb{F}I + X\otimes
X_-^{\bot})=\mathbb{F} I+(0)_+\otimes Y^*$ and there exists a ring
automorphism $\tau: \mathbb{F} \rightarrow \mathbb{F}$ and a map
$\gamma: X\otimes X_-^{\bot}\rightarrow \mathbb{F}$ such that
$$\Phi(x\otimes f) =\gamma(x,f) I + Df\otimes Cx \eqno(5.14)$$ holds for
all $x\in X$ and $f\in X_-^{\bot}$, where   $C: X\rightarrow Y^*$
and $D: X_-^{\bot} \rightarrow (0)_+$ are two $\tau$-linear
bijective maps. By a similar argument to the above,  one can check
that there exists an additive functional $h$ and a ring
anti-isomorphism $\Psi$ such that $\Phi(A)=-\Psi(A) + h(A)I$ for all
$A\in {\rm
Alg}\mathcal{N}$. 

This completes the proof of Theorem 1.1 for the case that
$(0)<X_-<X$ and $X_-$ is not complemented.\hfill$\Box$

The proof of Theorem 1.1 for the case that $(0) < (0)_+<X$ and
$(0)_+$ is not complemented  is similar and we omit it here.

Now, combining Sections 3-5, the proof of Theorem 1.1 is finished.

{\bf Acknowledgement.} The authors  wish to express their thanks to
the referees for   many comments to improve    the original draft.

\end{document}